\documentclass[10pt]{article}
\usepackage{amsthm}
\usepackage{latexsym}
\usepackage{amsmath}
\usepackage{amssymb}
\usepackage{epsfig}
\usepackage[matrix,arrow,cmtip,curve,dvips]{xy}

    \newtheorem{thm}{Theorem}[section]
    \newtheorem{prop}[thm]{Proposition}
    \newtheorem{lemma}[thm]{Lemma}

    \newtheorem{defn}[thm]{Definition}

    \newtheorem{rem}[thm]{Remark}
    \newtheorem{example}[thm]{Example}
    \newtheorem{question}[thm]{Question}

\newcommand{\id}{\mbox{\upshape id}}

\newcommand{\Hom}{\mbox{\upshape Hom}}
\newcommand{\conn}{\mbox{\upshape conn}}

\newcommand{\Map}{\mbox{\upshape Map}}
\newcommand{\im}{\mbox{\upshape im}}
\newcommand{\diam}{\mbox{\upshape diam}}

\title{Hom complexes and homotopy theory \\
in the category of graphs}
\date{June 10, 2008}

\author{Anton Dochtermann \thanks{Research supported in part by NSF grant
DMS-9983797}\\
\small Institut f\"{u}r Mathematik, MA 6-2 \\[-0.8ex]
\small Technische Universit\"{a}t Berlin \\[-0.8ex]
\small Strasse des 17. Juni 136 \\[-0.8ex]
\small 10623 Berlin, Germany \\[-0.8ex]
\small \texttt{anton.dochtermann@gmail.com} }

\begin{document}
\maketitle

\begin{abstract}
We investigate a notion of $\times$-homotopy of graph maps that is
based on the internal hom associated to the categorical product in
the category of graphs.  It is shown that graph $\times$-homotopy is
characterized by the topological properties of the $\Hom$ complex, a functorial way to assign a poset (and hence topological space) to a
pair of graphs; $\Hom$ complexes were introduced by Lov\'{a}sz and further studied by Babson and Kozlov to give topological bounds on chromatic number.  Along the way, we also establish some structural
properties of $\Hom$ complexes involving products and exponentials
of graphs, as well as a symmetry result which can be used to reprove
a theorem of Kozlov involving foldings of graphs.  Graph
$\times$-homotopy naturally leads to a notion of homotopy
equivalence which we show has several equivalent characterizations.
We apply the notions of $\times$-homotopy equivalence to the class
of dismantlable graphs to get a list of conditions that again
characterize these. We end with a discussion of graph homotopies
arising from other internal homs, including the construction of
`$A$-theory' associated to the cartesian product in the category of
reflexive graphs.
\end{abstract}

\section{Introduction}
In many categories, the notion of a pair of homotopic maps can be
phrased in terms of a map from some specified object into an
\textit{exponential object} associated to an internal hom structure
on that category (we will review these constructions below).  The
typical example is the category of (compactly generated) topological
spaces, where a homotopy between maps $f:X \rightarrow Y$ and $g:X
\rightarrow Y$ is nothing more than a map from the interval $I$ into
the topological space $\Map(X,Y)$. Other examples include simplicial
objects, as well as the category of chain complexes of $R$-modules.
For the latter, a chain homotopy between chain maps $f:C \rightarrow
D$ and $g:C \rightarrow D$ can be recovered as a map from the chain
complex $I$ (defined to be the complex consisting of 0 in all
dimensions except $R$ in dimensions 0 and 1, with the identity map
between them) into the complex $\Hom(C,D)$.

In this paper we consider these constructions in the context of the
category of graphs.  In particular, we investigate a notion of what
we call $\times$-homotopy that arises from consideration of the well
known internal hom associated to the categorical product.  Here the relevant construction is the exponential $H^G$, a graph whose looped vertices parameterize the graph homomorphisms (maps) from $G$ to $H$.  We use the
notion of (graph theoretic) connectivity to provide a notion of a
`path' in the exponential graph.  It turns out that
$\times$-homotopy classes of maps are related to the topology of the
so-called $\Hom$-complex, a functorial way to assign a poset
$\Hom(G,H)$ (and hence topological space) to a pair of graphs $G$
and $H$.  $\Hom$ complexes were first introduced by
Lov\`{a}sz in his celebrated proof of the Kneser conjecture (see
\cite{Lov78}), and were later developed by Babson and Kozlov in
their proof of the Lov\`{a}sz conjecture (see \cite{BKcom} and
\cite{BKpro}).  Elements of the poset $\Hom(G,H)$ are graph
\textit{multi}-homomorphisms from $G$ to $H$, with the set of graph
homomorphisms the atoms.  Fixing one of the two coordinates of the
$\Hom$ complex in each case provides a functor from graphs to
topological spaces, and in Theorem \ref{homotopic} we show that $\times$-homotopy of graph maps
can be characterized by the topological homotopy type of the maps
induced by these functors.

Graph $\times$-homotopy of maps naturally leads us to a notion of
homotopy equivalence of graphs, which in Theorem \ref{equivalent} we show can again be
characterized in terms of the topological properties of relevant
$\Hom$ complexes.  This result also exhibits a certain symmetry in the two entries of the $\Hom$ complex and can be used to reprove a result of Kozlov from \cite{Kfold}, here stated as Proposition \ref{folds}.  The graph operations known as `folding' and
`unfolding' preserve homotopy type, and in fact we show that in some
sense these operations generate the homotopy equivalence class of a
given graph.  In particular, a pair of \textit{stiff} graphs are
homotopy equivalent if and only if they are isomorphic.  One
particular case of interest arises when the graph can be folded down
to a single looped vertex, a class of graphs called
\textit{dismantlable} in the literature (see for example \cite{HN90}).  We apply Theorem \ref{equivalent} to obtain several characterizations of dismantlable graphs which adds to the list established by Brightwell and Winkler in \cite{BW00}.

The paper is organized as follows.  In Section 2 we describe the
category of graphs, and gather together some facts regarding its
structure. Here we focus on the internal hom structure associated
with the categorical product, and review the construction of the
exponential graph $H^G$ that serves as the right adjoint.  In
Section 3 we recall the construction of the $\Hom$ complex and discuss some properties.  We establish some structural results regarding preservation of homotopy type of the $\Hom$ complex under graph exponentiation (Proposition \ref{adjoint}) as well as arbitrary limits (e.g. products) of graphs (Proposition \ref{product}).  The latter has applications to special cases of Hedetniemi's conjecture, while the former allows us to interpret the complex $\Hom(G,H)$ in terms of the clique complex of
the exponential graph $H^G$.  It is this characterization that will
allow us to relate the topology of the $\Hom$ complex with
$\times$-homotopy classes of graph maps in later sections.

In Section 4 we introduce the notion of $\times$-homotopy of graph
maps in terms of paths in the exponential graph and work out some
examples.  We discuss the characterization of $\times$-homotopy in
terms of the topology of the relevant $\Hom$ complex.  The
construction of $\times$-homotopy naturally leads us to a graph
theoretic notion of \textit{homotopy equivalence} of graphs, and in
Section 5 we prove some equivalent characterizations in terms of the
topology of the $\Hom$ complexes.  In particular, it is the symmetry involved in this characterization that allows us to reprove the result of Kozlov discussed above.  We also discuss some of the
categorical properties that are satisfied.  In Section 6, we
investigate some of the structure of these homotopy equivalence
classes, and discuss the relationship with the graph operations
known as foldings and unfoldings and the related notion of a stiff
graph.  Here we apply our previous results to obtain several characterizations of the class of dismantlable graphs.  Finally, in Section 6, we
briefly discuss one other notion of homotopy that arises from the
internal hom associated to the cartesian product.  It turns out that
this construction recovers the existing notion of the so-called
$A$-theory of graphs discussed in \cite{BL}.

\noindent \textbf{Acknowledgments:} The author would like to thank
Matt Kahle for many helpful conversations, and an anonymous referee for suggestions.  Special thanks to his
advisor, Eric Babson, for his many insights into the questions
addressed in this paper.

\section{The category of graphs}

We will work in the category of graphs.  A \textit{graph} $G =
(V(G), E(G))$ consists of a vertex set $V(G)$ and an edge set $E(G)
\subseteq V(G) \times V(G)$ such that if $(v,w) \in E(G)$ then
$(w,v) \in E(G)$.  Hence our graphs are undirected and do not have
multiple edges, but may have loops (if $(v,v) \in E(G)$).  If $(v,w)
\in E(G)$ we will say that $v$ and $w$ are \textit{adjacent} and
denote this as $v \sim w$.  Given a pair of graphs $G$ and $H$, a
\textit{graph homomorphism} (or \textit{graph map}) is a mapping of
the vertex set $f:V(G) \rightarrow V(H)$ that preserves adjacency:
if $v \sim w$ in $G$, then $f(v) \sim f(w)$ in $H$ (equivalently
$(v,w) \in E(G)$ implies $(f(v),f(w)) \in E(H)$).  With these as our
objects and morphisms we obtain a category of graphs which we will
denote ${\mathcal G}$.  If $G$ and $H$ are graphs, we will use
${\mathcal G}(G,H)$ to denote the set of graph maps between them.

In this section, we review some of the structure of ${\mathcal G}$.
Of particular importance for us will be the existence of an
\textit{internal hom} associated to the categorical product.  We
start by recalling some related constructions, all of which can be
found in \cite{GR01} and \cite{HN04}.  For
undefined categorical terms, we refer to \cite{Mac98}.

\begin{defn}
For graphs $G$ and $H$, the \textit{categorical coproduct} $G \amalg
H$ is the graph with vertex set $V(G) \amalg V(H)$ and with
adjacency given by $(x,x^\prime) \in E(G \amalg H)$ if $(x,x^\prime)
\in E(G)$ or $(x,x^\prime) \in E(H)$.
\end{defn}

\begin{defn}
For graphs $G$ and $H$, the \textit{categorical product} $G \times
H$ is a graph with vertex set $V(G) \times V(H)$ and adjacency given
by $(g,h) \sim (g^\prime, h^\prime)$ in $G \times H$ if both $g \sim
g^\prime$ in $G$ and $h \sim h^\prime$ in $H$ (see Figure 1).
\end{defn}

\begin{center}
\epsfig{file=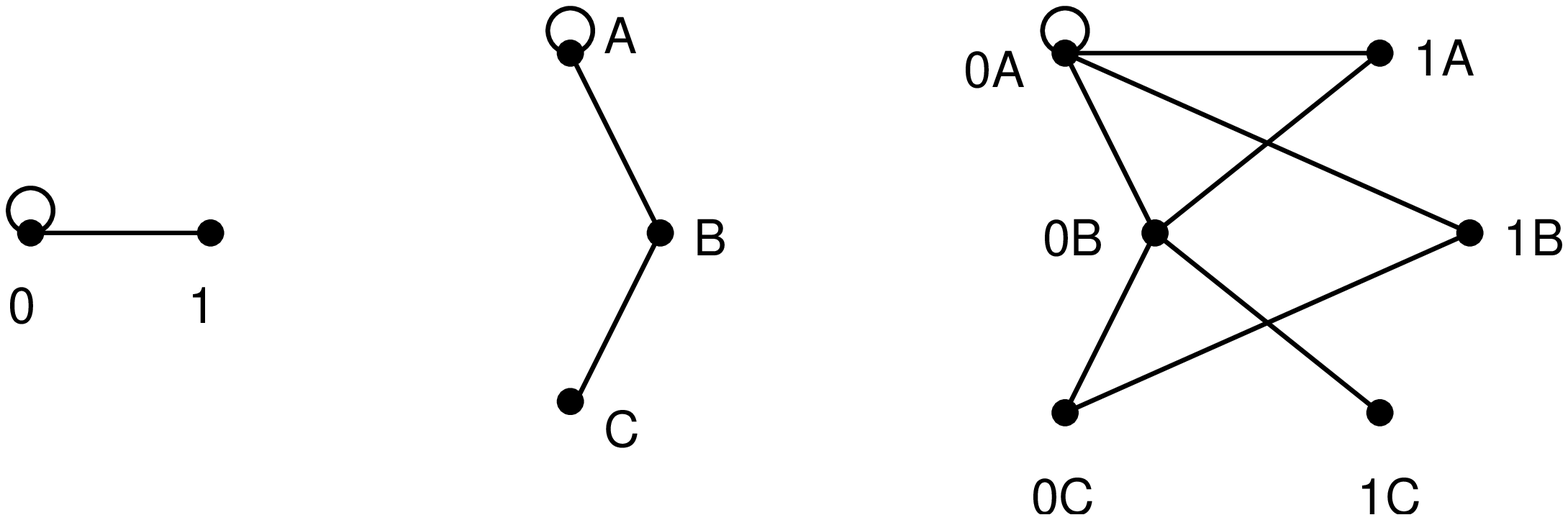, height=1.2 in, width = 3.5 in}

{Figure 1: The graphs $G$, $H$, and $G \times H$}

\end{center}

\begin{defn}
For graphs $G$ and $H$, the \textit{categorical exponential
graph} $H^G$ is a graph with vertex set $\{f:V(G) \rightarrow
V(H)\}$, the collection of all vertex set maps, with adjacency given
by $f \sim f^\prime$ if whenever $v \sim v^\prime$ in $G$ we have
$f(v) \sim f^\prime(v^\prime)$ in $H$ (see Figure 2).
\end{defn}

\vspace{.1 in}

\begin{center}
\epsfig{file=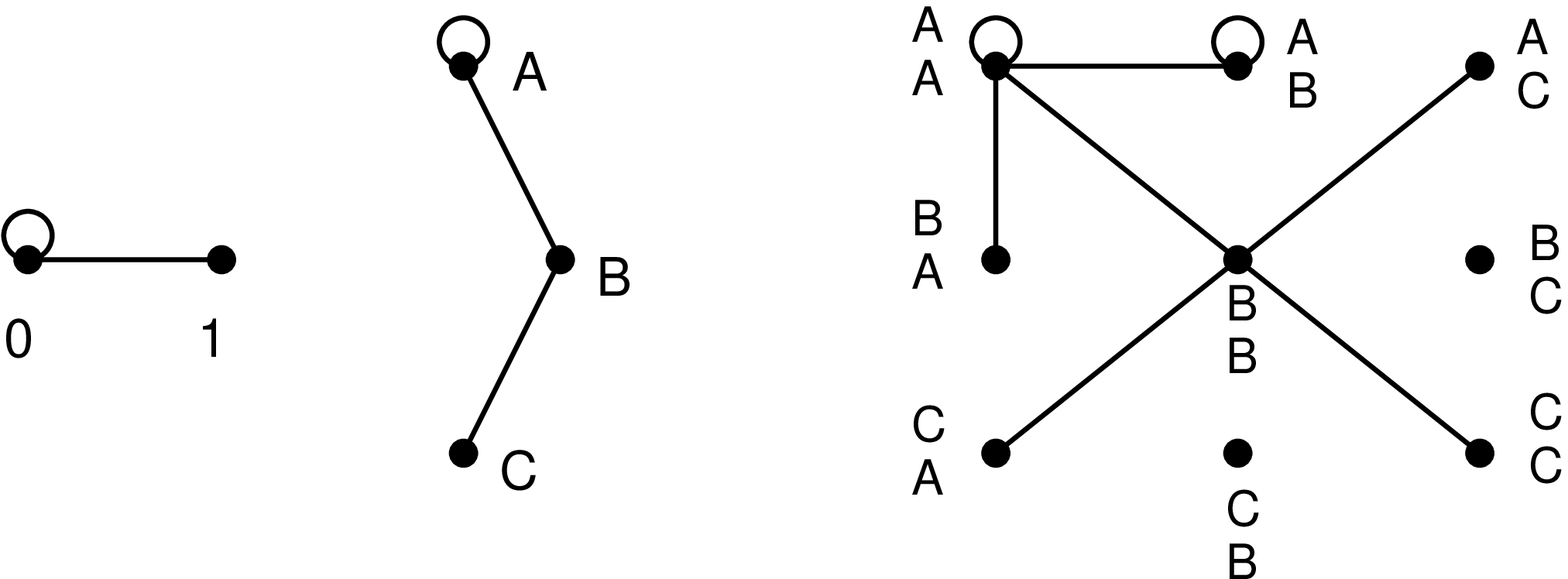, height=1.3 in, width = 3.5in}

{Figure 2: The graphs $G$, $H$, and $H^G$}
\end{center}

The next lemma shows that the exponential graph construction
provides a right adjoint to the categorical product.  By definition,
this gives the category of graphs the structure of an
\textit{internal hom} associated with the (monoidal) categorical
product.  This result is well known, and is more or less contained
in \cite{GR01}, but we state it here in a way that is consistent
with our notation.

\begin{lemma} \label{internal}
For graphs $A, B$ and $C$, we have a natural isomorphism of sets
\begin{displaymath}
\varphi:{\mathcal G}(A \times B, C) \rightarrow {\mathcal G}(A, C^B)
\end{displaymath}
given by $(\varphi(f)(v))(w) = f(v,w)$ for all $f \in {\mathcal G}(A
\times B, C)$, $v \in V(A)$, $w \in V(B)$.
\end{lemma}

\begin{proof}
Let $f:A \times B \rightarrow C$ be an element of ${\mathcal G}(A
\times B, C)$.  To see that $\varphi(f) \in {\mathcal G}(A, C^B)$,
suppose that $a \sim a^\prime$ are adjacent vertices in $A$.
 We need $\varphi(f)(a)$ and $\varphi(f)(a^\prime)$ to be adjacent
vertices in $C^B$.  To check this, suppose $b \sim b^\prime$ in $B$.
 Then we have $\varphi(f)(a)(b) = f(a,b)$ and
$\varphi(f)(a^\prime)(b^\prime) = f(a^\prime, b^\prime)$, which are
adjacent vertices of $C$ since $f$ is a graph map.

To check naturality, suppose $f:A \rightarrow A^\prime$ and $g:C
\rightarrow C^\prime$ are graph maps.  We need to verify that the
following diagram commutes:

\begin{center}

$\xymatrix{ {\mathcal G}(A \times B, C) \ar[d]_{\varphi}  &
{\mathcal G}(A^\prime \times B, C^\prime) \ar[l]_{(f
\times B, g)} \ar[d]^{\varphi} \\
{\mathcal G}(A,C^B) & {\mathcal G}(A^\prime, (C^{\prime})^B)
\ar[l]^{(f, g^B)} }$

\end{center}

For this, let $\alpha \in {\mathcal G}(A^\prime \times B, C)$. Then
on the one hand we have $(\varphi(f \times B, g))(\alpha)(a)(b) = (f
\times B, g)(\alpha)(a,b) = g(\alpha(f(a),b))$.  In the other
direction, we have $((f,g^B)(\varphi))(\alpha)(a)(b) =
g(\varphi(\alpha)(f(a))(b)) = g(\alpha(f(a),b))$.  Hence the diagram
commutes, and so the isomorphism $\varphi$ is natural.
\end{proof}

We close this section with a few additional definitions and remarks.  We let ${\bf 1}$ denote the graph consisting of a single looped vertex.  We point out that ${\bf 1}$ is the \textit{terminal object} in ${\mathcal G}$ in the sense that there exists a unique map $G \rightarrow {\bf 1}$ for all $G$.  Similarly, the graph $\emptyset$ is the graph whose vertex set is the empty set.  It is the \textit{initial object} in the sense that there exists a unique map $\emptyset \rightarrow G$ for all $G$.

A \textit{reflexive} graph $G$ is a graph with loops on all its
vertices ($v \sim v$ for all $v \in V(G)$).  A map of reflexive
graphs will be a graph map on the underlying graph.  We will use
${\mathcal G}^\circ$ to denote the category of reflexive graphs.

We see that ${\mathcal G}^\circ$ is a subcategory of ${\mathcal G}$,
and we let $i:{\mathcal G}^\circ \rightarrow {\mathcal G}$ denote
the inclusion functor.  Let $S:{\mathcal G} \rightarrow {\mathcal
G}^\circ$ denote the functor given by taking the subgraph induced by
looped vertices, and $L:{\mathcal G} \rightarrow {\mathcal G}^\circ$
denote the functor given by adding loops to all vertices (see Figure
3). One can check that $i$ is a left adjoint to $S$, whereas $i$ is
a right adjoint to $L$.  As functors ${\mathcal G} \rightarrow
{\mathcal G}$, one can check that $L$ (strictly speaking $iL$) is a
left adjoint to $S$ (strictly speaking $iS$).  We will make some use
of these facts in a later section.

\vspace{.1 in}

\begin{center}
\epsfig{file=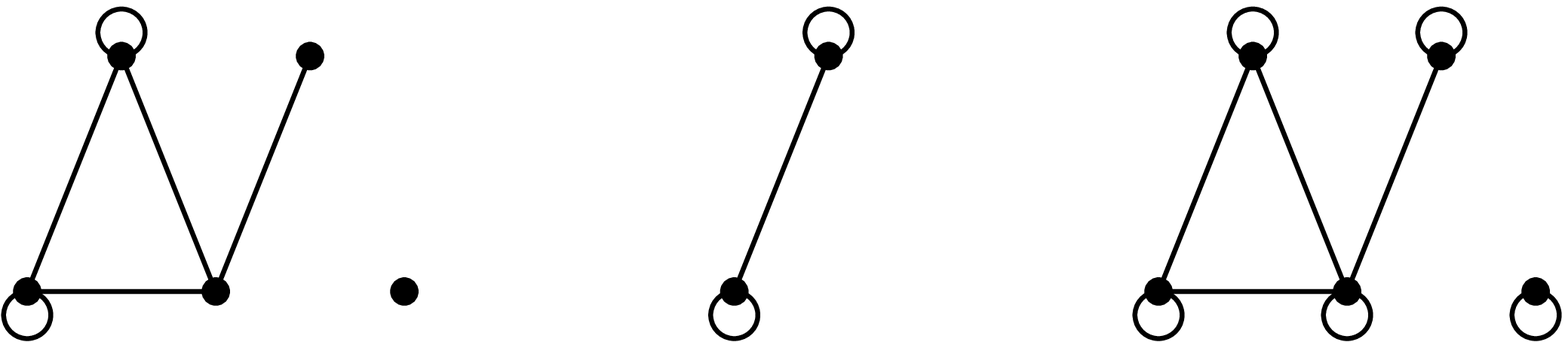, height=.7 in, width = 4 in}

{Figure 3: The graph $G$, and the reflexive graphs $S(G)$ and
$L(G)$.}
\end{center}

If $v$ and $w$ are vertices of a graph $G$, the \textit{distance}
$d(v,w)$ is the length of the shortest path in $G$ from $v$ to $w$.
The $diameter$ of a finite connected graph $G$, denoted $\diam(G)$
is the maximum distance between two vertices of $G$. The
\textit{neighborhood} of a vertex $v$, denoted $N_G(v)$ (or $N(v)$
if the context is clear), is the set of vertices that are adjacent
to $v$ (so that $v \in N(v)$ if and only if $v$ has a loop).

There are several simplicial complexes one can associate with a
given graph $G$.  One such construction is the \textit{clique
complex} $\Delta(G)$, a simplicial complex with vertices given by
all \textit{looped} vertices of $G$, and with faces given by all
cliques (complete subgraphs) on the looped vertices of $G$.

\section{The Hom complex and some properties}

Next we recall the construction of the $\Hom$ complex associated to
a pair of graphs.  As discussed in the introduction, (a version of)
the $\Hom$ complex was first introduced by Lov\'{a}sz in
\cite{Lov78}, and later studied by Babson and Kozlov in \cite{BKcom}.

\begin{defn}
For graphs $G, H$, we define $\Hom(G,H)$ to be the poset whose
elements are given by all functions $\eta: V(G) \rightarrow 2^{V(H)}
\backslash \{\emptyset\}$, such that if $(x,y) \in E(G)$ then
$(\tilde x, \tilde y) \in E(H)$ for every $\tilde x \in \eta(x)$ and
$\tilde y \in \eta(y)$.  The relation is given by containment, so
that $\eta \leq \eta^\prime$ if $\eta(x) \subseteq \eta^\prime(x)$
for all $x \in V$ (see Figure 4).
\end{defn}

\begin{center}
\epsfig{file=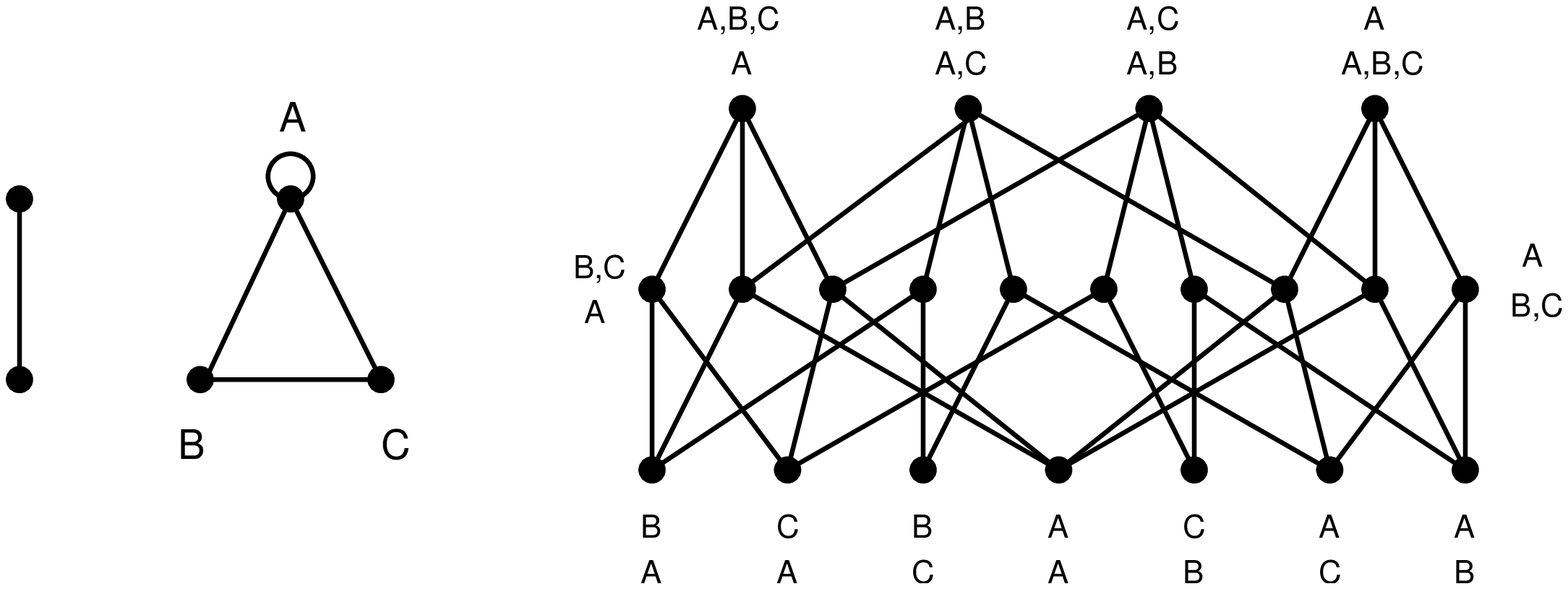, height=1.7 in, width = 5 in}

{Figure 4: The graphs $G$ and $H$, and the poset $\Hom(G,H)$}

\end{center}

We will often refer to $\Hom(G,H)$ as a topological space; by this
we mean the geometric realization of the order complex of the
poset.  The \textit{order complex} of a poset $P$ is the simplicial complex whose faces are the chains of $P$ (see Figure 5).

\begin{center}
\epsfig{file=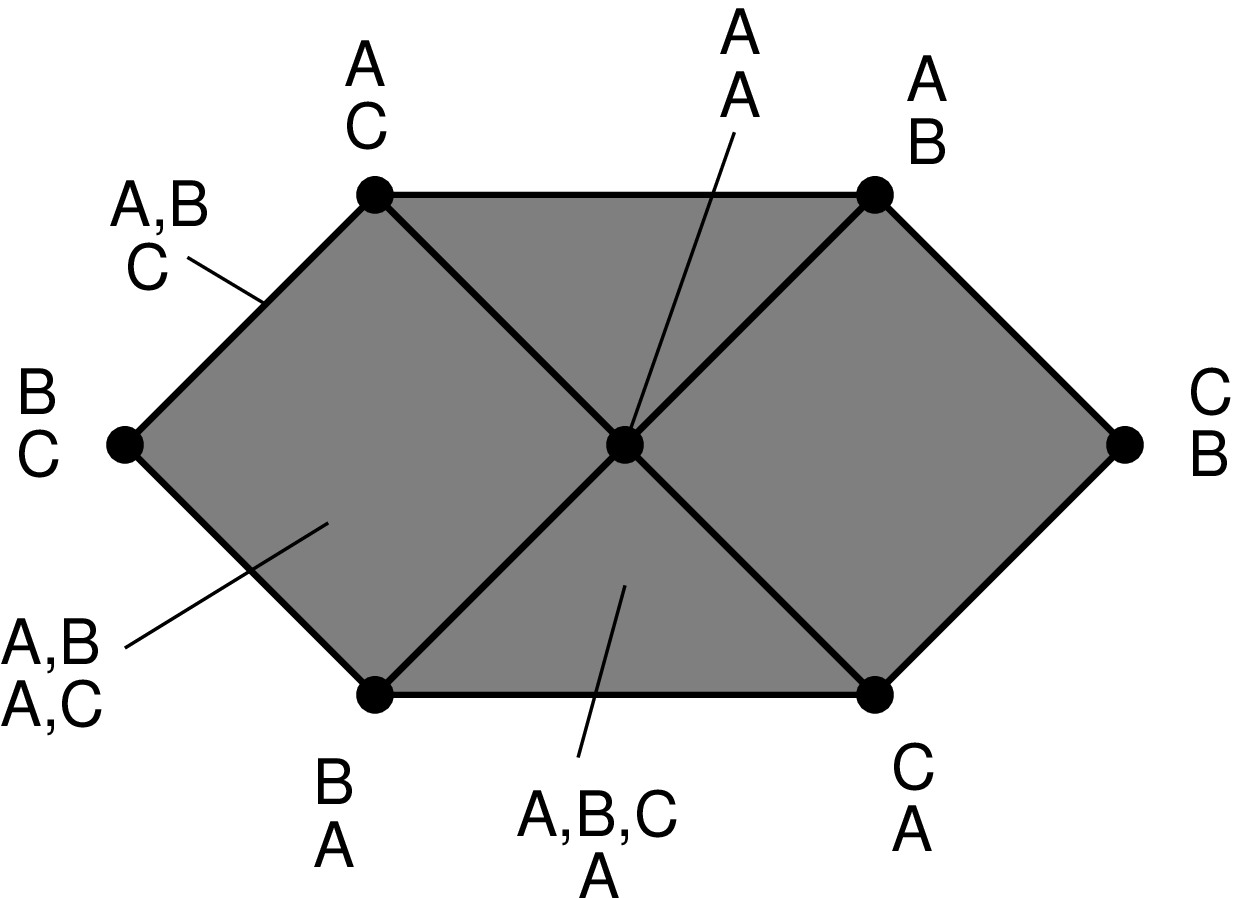, height=1.5 in, width = 1.75 in}

{Figure 5: The realization of the poset $\Hom(G,H)$ (up to
barycentric subdivision)}

\end{center}

Note that if $G$ and $H$ are both finite, then (the order complex of) this $\Hom(G,H)$ yields a
simplicial complex which is isomorphic to the barycentric subdivision
of the polyhedral $\Hom$ complex as defined in \cite{BKcom}.

The $\Hom$ complexes were originally used to obtain `topological' lower bounds on the chromatic number of graphs.  The main results of \cite{Lov78} and \cite{BKpro} in this context are the following theorems.  Here $\chi(G)$ is the chromatic number of a graph $G$, $\conn(X)$ denotes the (topological) connectivity of the space $X$, and $C_{2r+1}$ is the odd cycle of length $2r+1$.

\begin{thm}[Lov\'{a}sz] \label{LovTheorem}
For any graph $G$,
\[\chi(G) \geq \conn \big(\Hom(K_2, G)\big) + 3.\]
\end{thm}

\begin{thm}[Babson and Kozlov]
For any graph $G$,
\[\chi(G) \geq \conn \big(\Hom(C_{2r+1}, G)\big) + 4.\]
\end{thm}

In \cite{BKcom} Babson and Kozlov establish some basic functorial
properties of the $\Hom$ complex which we briefly discuss. Fixing
one of the coordinates of the $\Hom$ complexes provides a covariant
functor $\Hom(T,?)$, and a contravariant functor $\Hom(?,T)$, from
${\mathcal G}$ to the category of posets. If $f:G \rightarrow H$ is
a graph map, we have in the first case an induced poset map
$f_T:\Hom(T,G) \rightarrow \Hom(T,H)$ given by $f_T(\alpha)(t) =
\{f(g):g \in \alpha(t)\}$ for $\alpha \in \Hom(T,G)$ and $t \in
V(T)$. In the other case, we have $f^T:\Hom(H,T) \rightarrow
\Hom(G,T)$ given by $f^T(\beta)(g) = \beta(f(g))$ for $\beta \in
\Hom(H,T)$ and $g \in V(G)$.

A graph map $f: G \rightarrow H$ induces a natural transformation
$\bar f: \Hom(?,G) \rightarrow \Hom(?,H)$ in the following way.  For
each $T \in$ Ob(${\mathcal G}$) we have a map $\bar f_T : \Hom(T,G)
\rightarrow \Hom(T,H)$ given by $(\bar f_T (\alpha))(t) = \{f(g): g
\in \alpha(t)\}$ for $\alpha \in \Hom(T,G)$ and $t \in V(T)$. If
$g:S \rightarrow T$ is a graph map, the diagram

\begin{center}
$\xymatrix{ \Hom(S,G) \ar[d]^{\bar f_S}
& \Hom(T,G) \ar[l]^{g^G} \ar[d]^{\bar f_T} \\
\Hom(S,H) & \Hom(T,H) \ar[l]^{g^H} }$
\end{center}

\noindent commutes since if $\alpha \in \Hom(T,G)$ and $s \in V(S)$
then on the one hand we have
\[((\bar f_S g^G)(\alpha))(s) = \{f(x):x \in ((g^G
(\alpha))(s) \} = \{f(x): x \in \alpha(g(s))\},\] \noindent and on
the other
\[((g^H \bar f_T)(\alpha))(s) = ((\bar f_T)(\alpha))(g(s))
= \{f(x): x \in \alpha(g(s))\}.\]

The function induced by composition $\Hom(G,H) \times \Hom(H,K)
\rightarrow \Hom(G,K)$ is a poset map; see \cite{Kchr} for a proof
of this fact.

Many operations in the category of graphs interact nicely with the
topology of the $\Hom$ complexes.  We now gather together some of
these results.  The first observation comes from \cite{BKcom}.

\begin{lemma}
Let $A$, $B$, and $C$ be graphs.  Then there is an isomorphism of
posets
\begin{displaymath}
\Hom(A \amalg B, C) \cong \Hom(A,C) \times \Hom(B,C).
\end{displaymath}
Also, if $A$ is connected and not a single vertex, then
\begin{displaymath}
\Hom(A, B \amalg C) \cong \Hom(A,B) \amalg \Hom(A,C).
\end{displaymath}
\end{lemma}

As we will see, other graph operations are preserved by the $\Hom$
complexes \textit{up to homotopy type}.

Recall that for graphs $A, B$, and $C$ the exponential graph
construction provides the adjunction ${\mathcal G}(A \times B, C) =
{\mathcal G}(A, C^B)$, an isomorphism of sets.  The next proposition
shows that this map induces a homotopy equivalence of the associated
$\Hom$ complexes.  In the proof of Propositions \ref{adjoint} and \ref{product} we will use the following notion from poset topology (see \cite{Bjo95} for a good reference).  If $P$ is a poset, and $c:P \rightarrow P$ is a poset map such that $c \circ c = c$ and $c(p) \geq p$ for all $p \in P$ then $c$ is called a \textit{closure map}.  It can be shown (see \cite{Bjo95}) that in this case $c:P \rightarrow c(P)$ induces a strong deformation retract of the associated spaces.

\begin{prop} \label{adjoint}
Let $A$, $B$, and $C$ be graphs.  Then $\Hom(A \times B, C)$ can be
included in $\Hom(A, C^B)$ so that $\Hom(A \times B, C)$ is the
image of a closure map on $\Hom(A, C^B)$.  In particular, there is
an inclusion of a strong deformation retract

\begin{center}

$\xymatrix{ |\Hom(A \times B, C)| \ar@{^{(}->}[r]_{\hspace{.1
in}\simeq} & |\Hom(A, C^B)|. }$

\end{center}

\end{prop}

\begin{proof}
Let $P = \Hom(A \times B, C)$ and $Q = \Hom(A, C^B)$ be the
respective posets.  Our plan is to define an inclusion map $j:P
\rightarrow Q$ and a closure map $c:Q \rightarrow Q$ such that
$\im(j) = \im(c)$, from which the result would follow.

We define a map of posets $j: P \rightarrow Q$ according to
\[j(\alpha)(a) = \big\{f:V(B)\rightarrow V(C) | f(b) \in \alpha(a,b)
~\forall b \in B \big\},\]
for every $a \in V(A)$ and $\alpha \in P$.  To show that
$j(\alpha)$ is in fact an element of $Q$, we need to verify that if
$a \sim a^\prime$ in $A$ then we have $f \sim f^\prime$ in $C^B$ for all
$f \in j(\alpha)(a), f^\prime \in j(\alpha)(a^\prime)$.  If $b \sim b^\prime$ in $B$ then
$(a,b) \sim (a^\prime,b^\prime)$ in $A \times B$.  Hence $c \sim c^\prime$ in $C$ for $c \in \alpha(a,b)$ and all $c^\prime
\in \alpha(a^\prime, b^\prime)$.  In particular, $f(b) \sim f(b^\prime)$ in $C$ and we conclude that $f \sim f^\prime$, as desired.  Hence $j(\alpha)$ is indeed an element of $\Hom(A, C^B)$.

We claim that $j$ is injective.  To see this, let $\alpha \neq
\alpha^\prime$ be distinct elements of the poset $P$, with
$\alpha(a,b) \neq \alpha^\prime(a,b)$ for some $(a,b) \in V(A \times
B)$.  Without loss of generality, suppose $c \in \alpha(a,b)
\backslash \alpha^\prime(a,b)$.  Then we have some $f \in j(\alpha)(a)$
such that $f(b) = c$, and yet $f \notin j(\alpha^\prime)(a)$.  We
conclude that $j(\alpha) \neq j(\alpha^\prime)$, and hence $j$ is
injective.

Next we define a closure map of posets $c: Q \rightarrow Q$.  If
$\gamma:V(A) \rightarrow 2^{V(C^B)} \backslash \{\emptyset \}$ is an
element of $Q = \Hom(A, C^B)$, define $c(\gamma) \in Q$ as follows:
fix some $a \in V(A)$, and for every $b \in V(B)$ let $C_{ab}^\gamma
= \{f(b) \in V(C): f \in \gamma(a)\}$.  Then define $c(\gamma)$
according to
\[c(\gamma)(a) = \big\{g:V(B) \rightarrow V(C) | g(b) \in C_{ab}^\gamma ~\forall b \in B \big\}.\]

We first verify that $c$ maps into $Q$, so that $c(\gamma) \in \Hom(A,
C^B)$.  For this suppose $a \sim a^\prime$ in $A$ and let $f \in c(\gamma)(a)$ and $g \in c(\gamma)(a^\prime)$.
To show that $f \sim g$ in $C^B$ we consider some $b \in b^\prime$ in $B$.  Then by construction there is some $f^\prime \in \gamma(a)$,
$g^\prime \in \gamma(a^\prime)$ such that $f(b) = f^\prime(b)$ and
$g(b^\prime) = g^\prime (b^\prime)$. Hence $f(b) \sim g(b^\prime)$ in $C$ as desired.

It is clear that $c(\gamma) \geq \gamma$ and $(c \circ c)(\gamma) = c(\gamma)$ for all
$\gamma \in Q$.  Thus $c$ is a closure map.

Next we claim that $c(Q) \subseteq j(P)$.  To see this, suppose $\gamma \in Q$.  We define $\alpha: V(A \times B)
\rightarrow 2^{V(C)} \backslash \{\emptyset \}$ according to $\alpha(a,b) =
C_{ab}^\gamma$, where $C_{ab}^\gamma \subseteq V(C)$ is as
above.  To see that $\alpha \in \Hom(A \times B, C)$ suppose $(a,b) \sim (a^\prime, b^\prime)$, and let $c \in
\alpha(a,b) = C_{ab}^\gamma$, $c^\prime \in \alpha(a^\prime, b^\prime)
= C_{a^\prime b^\prime}^\gamma$.  Since $a \sim a^\prime$ we have $f \sim f^\prime$ in $C^B$ for all $f \in \gamma(a), f^\prime
\in \gamma(a^\prime)$.  Hence since $b \sim b^\prime$ we get $f(b) \sim f^\prime(b^\prime)$, and in particular obtain
$c \sim c^\prime$ in $C$ as desired.

Finally, we get $j(P) \subseteq c(Q)$ since $j(P) \subseteq Q$ and $c(j(P)) =
j(P)$. Thus $j(P) = c(Q)$, implying that $\Hom(A \times B, C) \simeq
\Hom(A, C^B)$ via this inclusion.
\end{proof}

\begin{rem}
As a result of Proposition \ref{adjoint}, for all graphs $G$ and $H$ there is
a homotopy equivalence
\[ \Hom(G,H) = \Hom({\bf 1} \times G, H) \simeq \Hom({\bf 1}, H^G),\]
\noindent where {\bf 1} denotes a single looped vertex.  The last of
these posets is the face poset of the clique complex on the looped
vertices of $H^G$, and hence its realization is the barycentric
subdivision of the clique complex of $H^G$.  Since the looped vertices in
$H^G$ are precisely the graph homomorphisms $G \rightarrow H$, we
see that $\Hom(G,H)$ can be realized up to homotopy type as the
clique complex of the subgraph of $H^G$ induced by the graph
homomorphisms.
\end{rem}

By a \textit{diagram of graphs} $D = \{D_i\}$, we mean a collection
of graphs $\{D_i\}$ with a specified collection of maps between them
(the image of a category $D$ under some functor to ${\mathcal G}$).
For a graph $T$, any such diagram of graphs gives rise to a diagram
of posets $\Hom(T,D)$ obtained by applying the functor $\Hom(T,?)$
to each object and each morphism (see Figure 6).

\begin{center}
$\xymatrix{
& \cdot \ar[d] \\
\cdot \ar[r] & \cdot }$ \hspace {.75 in}
$\xymatrix{
& D_1 \ar[d]^f \\
D_3 \ar[r]_g & D_2 }$  \hspace {0.75 in}
 $\xymatrix{
& \Hom(T,D_1) \ar[d]^{f_T} \\
\Hom(T,D_3) \ar[r]_{g_T} & \Hom(T,D_2) }$

{Figure 6: A category $D$, a diagram of graphs, and the induced
diagram of posets}

\end{center}

We can combine the facts from Lemma \ref{internal} and Proposition
\ref{adjoint} to see that $\Hom$ complexes preserve (up to homotopy
type) limits of such diagrams.

\begin{prop} \label{limit}
Let $D$ be a diagram of graphs with limit $\lim (D)$. Then for
every graph $T$ we have a homotopy equivalence:
\begin{displaymath}
\big|\Hom\big(T,\lim (D)\big)\big| \simeq \big|\lim \big(\Hom(T,
D)\big)\big|.
\end{displaymath}
\end{prop}

\begin{proof}
Let $T$ be a graph.  We will express the functor $\Hom(T,?)$ as a
composition of functors that each preserve limits.  First we note
that the functor $(?)^T:{\mathcal G} \rightarrow {\mathcal G}$ given
by $G \mapsto G^T$ preserves limits since it has the left adjoint
given by the functor $? \times T$; this was the content of
Proposition \ref{adjoint}. Hence for any diagram of graphs $D$, we
get $\big(\lim (D)\big)^T = \lim (D^T)$.

Next we note that the functor $L:{\mathcal G} \rightarrow {\mathcal
G}^{\circ}$ that takes the induced subgraph on the looped vertices
(described above) also preserves limits since it has the left
adjoint given by the inclusion functor ${\mathcal G}^\circ
\rightarrow {\mathcal G}$. So we have $L\big(\lim (D)\big) = \lim \big(L(D)\big)$.

Now we claim that the functor $\Hom({\bf 1}, ?)$ preserves limits up
to homotopy type.  To see this, we recall that $\Hom({\bf 1}, ?)$,
as a functor from the category of reflexive graphs, associates to a
given (reflexive) graph $G$ the face poset of its \textit{clique
complex}, $\Delta(G)$. Hence, taking geometric realization, we get
$|\Hom({\bf 1}, G)| \simeq |\Delta(G)|$ for all reflexive graphs $G$. Now, as a functor to flag simplicial complexes, the clique
complex $\Delta$ has an inverse functor given by taking the
1-skeleton and adding loops to each vertex.  In particular, this
shows that $\Delta(?)$ preserves limits, and we get $\Delta(\lim
\tilde D) = \lim \big(\Delta(\tilde D)\big)$, for any diagram of reflexive
graphs $\tilde D$.

Finally, we can put these observations together to get the following
string of isomorphisms ($=$) and homotopy equivalences ($\simeq$):
\begin{align*}
\big|\Hom\big(T,\lim (D)\big)\big| &\simeq \big|\Hom\big({\bf 1},\big(\lim (D)\big)^T\big)\big| = \big|\Hom\big({\bf 1},\lim (D^T)\big)\big|\\
&= \Big|\Hom\Big({\bf 1}, L\big(\lim (D^T)\big)\Big)\Big| = \Big|\Hom\Big({\bf 1}, \lim\big(L(D^T)\big)\Big)\Big| \\
&\simeq \Big|\Delta\Big(\lim\big(L(D^T)\big)\Big)\Big| = \Big|\lim \Big(\Delta\big(L(D^T)\big)\Big)\Big| \\
&\simeq \Big| \lim\Big(\Hom \big({\bf 1}, L(D^T)\big)\Big)\Big| = \Big|\lim \Big(\big(\Hom \big({\bf 1}, (D^T)\big)\Big)\Big| \\
&\simeq \big| \lim \big(\Hom(T, D)\big)\big|.
\end{align*}

\noindent The first and last homotopy equivalences are as in
Proposition \ref{adjoint}.
\end{proof}

Recall that the product $G \times H$ is a limit (pullback) of the
diagram $G \rightarrow {\bf 1} \leftarrow H$.  Since $\Hom(T, {\bf
1})$ is a point for every graph $T$, this implies that $|\Hom(T,G)|
\times |\Hom(T,H)|$ is homotopy equivalent to $|\Hom(T,G \times
H)|$ for all graphs $T$, $G$, and $H$.  In fact in the case of the product we can exhibit this homotopy equivalence as a closure map on the level of posets.

\begin{prop} \label{product}
Let $T$, $G$, and $H$ be graphs.  Then the poset $\Hom(T,G) \times
\Hom(T,H)$ can be included into $\Hom(T,G \times H)$ so that
$\Hom(T,G) \times \Hom(T,H)$ is the image of a closure map on
$\Hom(T, G \times H)$. In particular, there is an inclusion of a
strong deformation retract

\begin{center}
$\xymatrix{ |\Hom(T,G)| \times |\Hom(T,H)|
\ar@{^{(}->}[r]_{\hspace{.35 in}\simeq} & |\Hom(T,G \times H)|. }$
\end{center}
\end{prop}

\begin{proof}
We let $Q = \Hom(T,G) \times \Hom(T,H)$ and $P = \Hom(T, G \times
H)$ be the respective posets.  Once again our plan is to define an
inclusion $i:Q \rightarrow P$ and a closure map $c:P \rightarrow P$
such that $\im(i) = \im(c)$.

We define a map $i:Q \rightarrow P$ according to $i(\alpha,
\beta)(v) = \alpha(v) \times \beta(v)$, for every vertex $v \in V(T)$.
Note that if $v$ and $w$ are adjacent vertices of $T$ then $\tilde v \sim \tilde w$ in $G$ and $v^\prime, w^\prime$ in $H$ for all
$\tilde v \in \alpha(v)$, $\tilde w \in \alpha (w)$, $v^\prime \in
\beta(v)$, and $w^\prime \in \beta(w)$.  Hence $(\tilde v, \tilde
w) \sim (v^\prime, w^\prime)$ are adjacent in $G \times H$, so that $i(\alpha,
\beta)$ is indeed an element of $\Hom(T, G \times H)$.  It is clear that $i$
is injective.

Next, we define a closure operator $c:P \rightarrow P$, whose
image will coincide with that of the map $i$.  For $\gamma \in P := \Hom(T, G \times H)$, we define $c(\gamma) \in P$ as follows: for every $v \in V(T)$ we have minimal vertex subsets $A_v \subseteq V(G)$, $B_v \subseteq V(H)$ such that $\gamma(v)
\subseteq \{(a, b): a \in A_v, b \in B_v\}$.  Define $c(\gamma)(v) :=
\{(a, b)\} = A_v \times B_v$ to be this minimal set of vertices of
$G \times H$.

We first verify that $c$ maps into $P$, so that $c(\gamma) \in \Hom(T,
G \times H)$.  Suppose $v \sim w$ are adjacent vertices of $T$.  If
$(\tilde a, \tilde b) \in c(\gamma)(v)$ and $(a ^\prime, b^\prime)
\in c(\gamma)(w)$ then we have $(\tilde a, \tilde y), (\tilde x,
\tilde b) \in \gamma(v)$ and $(a^\prime, y^\prime), (x^\prime,
b^\prime) \in \gamma(w)$ for some $\tilde x, x^\prime \in G$ and
$\tilde y, y^\prime \in H$.  Hence $\tilde a \sim a^\prime$ in $G$ and also $\tilde b, b^\prime$ in $H$, so that $(\tilde a, \tilde
b) \sim (a^\prime, b^\prime)$ in $G \times H$ as desired.

Since $c(\gamma) \geq \gamma$ and $(c \circ c)(\gamma) = c(\gamma)$ for all $\gamma \in P$,
we see that $c:P \rightarrow P$ is a closure operator.

Next we claim that $c(P) \subseteq i(Q)$.  Suppose $c(\gamma) \in
c(P)$, so that for all $v \in T$ we have $c(\gamma)(v) = A_v \times
B_v$ for some $A_v \subseteq V(G)$ and $B_v \subseteq V(H)$.  Define
$\alpha:V(T) \rightarrow  2^{V(G)} \backslash \{\emptyset\}$ by
$\alpha(v) = A_v$, and $\beta:V(T) \rightarrow 2^{V(H)} \backslash
\{\emptyset\}$ by $\beta(v) = B_v$.  We claim that $\alpha \in
\Hom(T,G)$ and $\beta \in \Hom(T,H)$.  Indeed, if $w \in T$ is a
vertex adjacent to $v$ and $\alpha(w) = A_w$, then if $a_i \in A_v$
and $a_{i^\prime} \in A_w$, we have $(a_i, y) \in \gamma(v)$ and
$(a_{i^\prime}, y^\prime) \in \gamma(w)$ for some $y, y^\prime \in
H$. Hence $(a_i, y)$ and $(a_{i^\prime}, y^\prime)$ are adjacent
vertices in $G \times H$ (since $\gamma \in \Hom(T,G \times H))$.
But this implies that $a_i$ and $a_{i^\prime}$ are adjacent in $G$,
as desired.

Finally, $i(Q) \subseteq c(P)$ since $i(Q) \subseteq P$ and $c(i(Q))
= i(Q)$.  Thus $i(Q) = c(P)$ and hence $\Hom(T,G) \times \Hom(T,H)
\simeq \Hom(T,G \times H)$  via this inclusion.
\end{proof}

\begin{rem}
Proposition \ref{product} can be used to prove special cases of \textit{Hedetniemi's conjecture}, which is the simple statement that $\chi(G \times H) = \min \{\chi(G), \chi(H)\}$ for all graphs $G$ and $H$.  Since it is clear that $\chi(G \times H) \leq \min \{\chi(G), \chi(H)\}$, the content of the conjecture is the other inequality.  Now, combining Proposition \ref{product} together with (say) Theorem \ref{LovTheorem} we obtain
\begin{align*}
\chi(G \times H) &\geq \conn \big(\Hom(K_2, G \times H)\big) + 3\\
&= \conn \big(\Hom(K_2,G) \times \Hom(K_2, H)\big) + 3\\
&= \min \big\{\conn \big(\Hom(K_2,G)\big), \conn \big(\Hom(K_2,H)\big) \big\} + 3.
\end{align*}
Here we apply the simple observation that $\conn(X \times Y) = \min\{\conn(X), \conn(Y)\}$ for topological spaces $X$ and $Y$.  This then proves the conjecture for the case when the topological bounds on the chromatic numbers of $G$ and $H$ are tight (e.g., when $G$ and $H$ are both taken to be either Kneser graphs or generalized Mycielski graphs).
\end{rem}

\section {Graph $\times$-homotopy and Hom complexes}

In this section, we define a notion of homotopy for graph maps and
describe its interaction with the $\Hom$ complexes.  The motivation
comes from the internal hom structure in the category ${\mathcal G}$
as described above.

Recall that a vertex set map $f:V(G) \rightarrow V(H)$ is a looped
vertex in $H^G$ if and only if $f$ is a \textit{graph} map $G
\rightarrow H$. Hence the set of graph maps ${\mathcal G}(G,H)$ are
precisely the looped vertices in the internal hom graph $H^G$. The
(path) connected components of the graph $H^G$ then provide a
natural notion of `homotopy' for graph maps: two maps $f,g:G
\rightarrow H$ will be considered $\times$-\textit{homotopic} if we
can find a path along the looped vertices $H^G$ that starts at $f$
and ends at $g$.  The use of the $\times$ is to emphasize the fact that we are using the exponential graph construction which is adjunct to the \textit{categorical product}; in the last section we will consider other exponentials.

To make the notion of a path truly graph theoretic we want to think
of it as a map from a path-like graph object into the graph $H^G$.

\begin{defn}
We let $I_n$ denote the graph with vertices $\{0,1, \dots, n\}$ and
with adjacency given by $i \sim i$ for all $i $ and $(i-1) \sim i$
for all $1 \leq i \leq n$ (see Figure 7).
\end{defn}

\vspace{.1 in}
\begin{center}
\epsfig{file=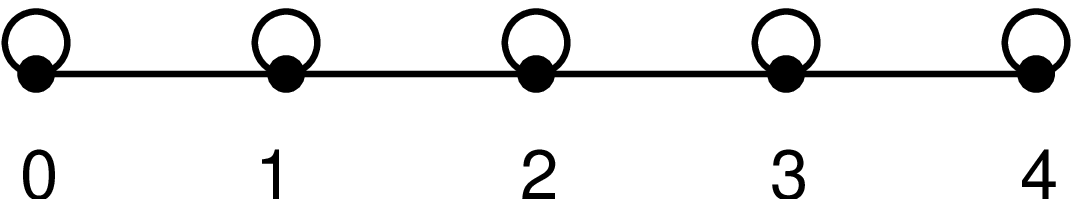, height=.3 in, width = 1.5 in}

{Figure 7: The graph $I_4$.}
\end{center}

Note that $N(n) = \{n, n-1\} \subseteq \{n, n-1, n-2\} = N(n-1)$, and
hence we can fold the endpoint of $I_n$.  This gives us the
following property.

\begin{lemma} \label{contractible}
$\Hom(T, I_n)$ is contractible for all $n \geq 0$ and every graph $T$.
\end{lemma}

\begin{proof}
We proceed by induction on $n$. For $n=0$, we have that $\Hom(T,
I_0) = \Hom(T, {\bf 1})$ is a point.  For $n > 0$, we use the fact
that $N(n) \subseteq N(n-1)$ to get $\Hom(T,I_n) \simeq \Hom(T,
I_{n-1})$ by Proposition \ref{folds}.  The latter complex is
contractible by induction.
\end{proof}

\begin{defn}
Let $f,g:G \rightarrow H$ be graph maps.  We say that $f$ and $g$
are $\times$-\textit{homotopic} if there exists an integer $n \geq
1$ and a map of graphs $F:I_n \rightarrow H^G$ such that $F(0) = f$
and $F(n) = g$.  In this case we will also say the maps are
\textit{$n$-homotopic}.
\end{defn}

We will denote $\times$-homotopic maps as $f \simeq_\times g$, or
simply $f \simeq g$ if the context is clear. Graph $\times$-homotopy
determines an equivalence relation on the set of graph maps between
$G$ and $H$, and we let $[G,H]_\times$ (or simply $[G,H]$) denote the set of
$\times$-homotopy classes of maps between graphs $G$ and $H$.

\begin{example}
As an example we can take $G = K_2$ and $H = K_3$ to be the complete
graphs on 2 and 3 vertices.  The graph $H^G$ is displayed in Figure
8.

\begin{center}
\epsfig{file=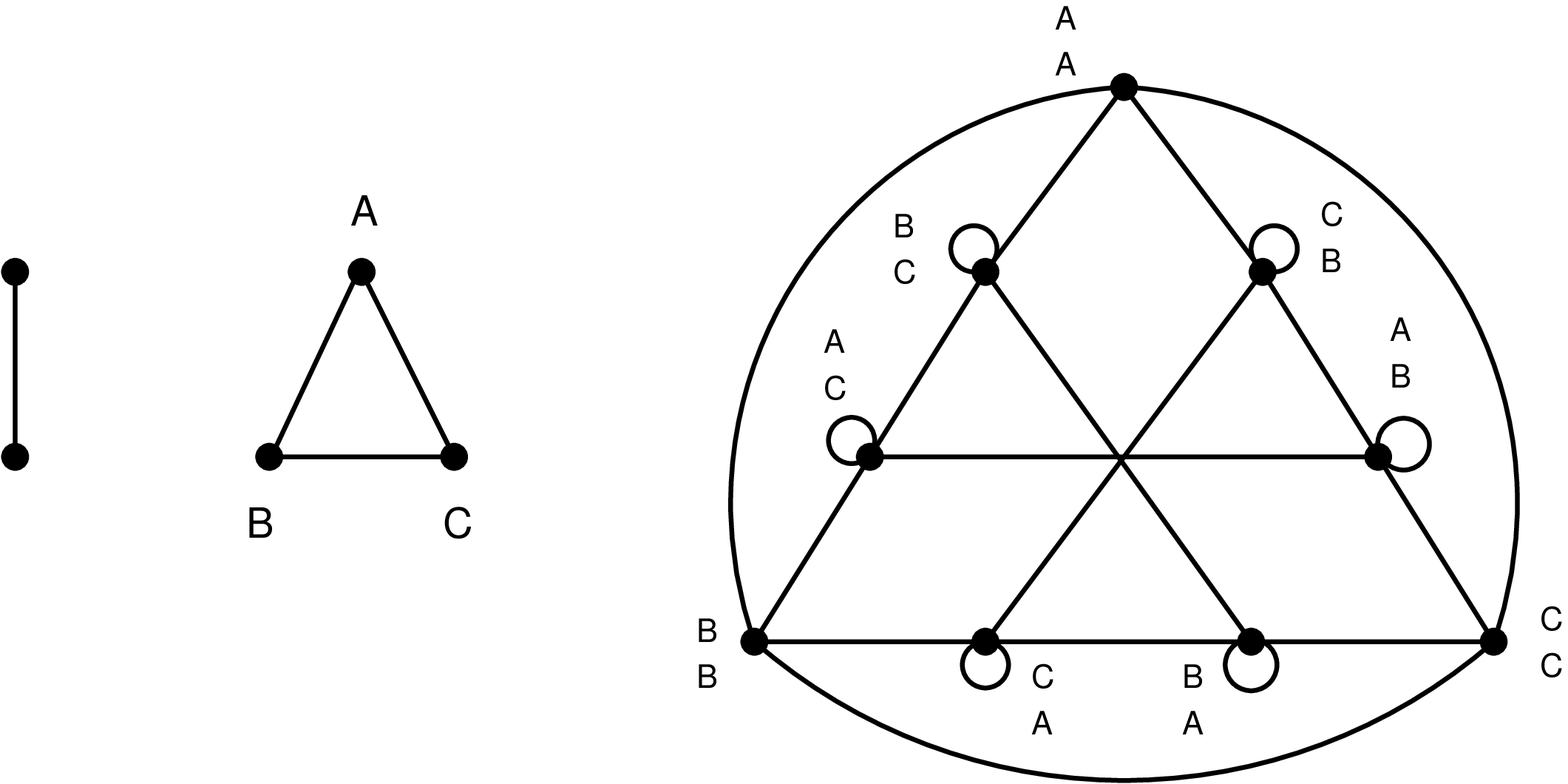, height = 2 in, width = 4.4 in}

{Figure 8: The graphs $G=K_2$, $H=K_3$, and $H^G$.}
\end{center}

We see that each of the six graph maps $f:G \rightarrow H$ is
represented by a looped vertex in the exponential graph $H^G$.  In
this case, any two maps $f$ and $g$ are connected by a path along
other looped vertices, and hence in our setup all maps from $G =
K_2$ to $H = K_3$ will be considered $\times$-homotopic (so that
there is a single homotopy class of maps).
\end{example}

\begin{example}
On the other hand, if we take $G=K_2$, and this time $H=K_2$, we get
two distinct $\times$-homotopy classes of maps.  The graph $H^G$ is
displayed in Figure 9.

\begin{center}
\epsfig{file=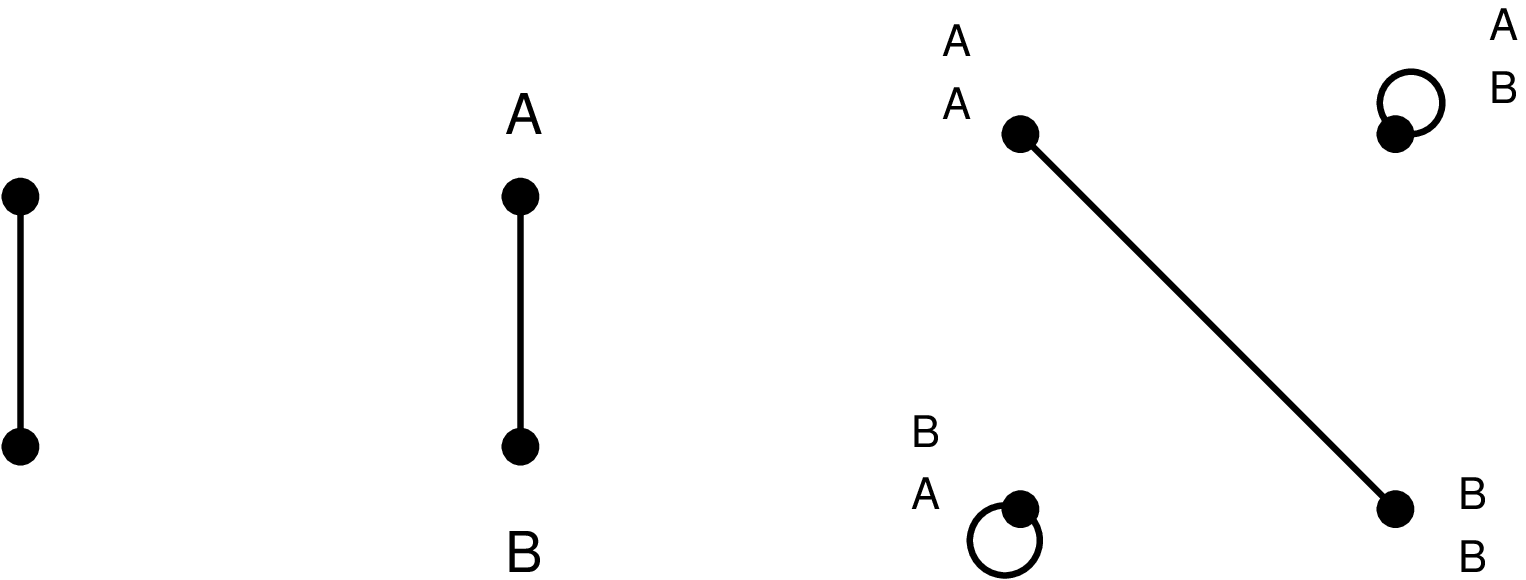, height = 1.1 in, width = 3 in}

{Figure 9: The graphs $G=K_2$, $H=K_2$, and $H^G$}

\end{center}

We see that the two graph maps $G \rightarrow H$ are represented by
looped vertices in $H^G$, but this time are disconnected from one
another.  Hence in this example, each of the two graph maps is in
its own $\times$-homotopy class.

\end{example}

We can understand $\times$-homotopy in other ways by considering the
adjoint properties available to us.  Note that for all $m \leq n$,
we have a map $\iota_m:G \rightarrow G \times I_n$ given by $v
\mapsto (v,m)$, an isomorphism onto its image.  A map $F:I_n
\rightarrow H^G$ corresponds to a map $\tilde F:G \times I_n
\rightarrow H$ with the property that $\tilde F \times 0 = f$ and
$\tilde F \times n = g$. It is this formulation that we will most
often use to check for $\times$-homotopy.  We record this
observation as a lemma.

\begin{lemma}
Let $f,g:G \rightarrow H$ be graph maps.  Then $f$ and $g$ are
$\times$-homotopic if and only if there exists an integer $n$ and a
graph map $F: G \times I_n \rightarrow H$ such that $F_0 := F
\circ \iota_0 = f: G \rightarrow H$ and $F_n := F \circ \iota_n
= g: G \rightarrow H$.

\begin{center}

$\xymatrix{
G \ar[d]_{\iota_0} \ar[dr]^{f} \\
G \times I_n \ar[r]^F & H \\
G \ar [u]^{\iota_n} \ar [ur]_{g} }$

\end{center}

\end{lemma}

Next we investigate how $\times$-homotopy of graph
maps interacts with the $\Hom$ complex.  It turns out that
$\times$-homotopy equivalence classes of maps are characterized by
the topology of the $\Hom$ complex in the following way.

\begin{prop} \label{component}
Let $G$ and $H$ be graphs, and suppose $f,g:G \rightarrow H$ are graph maps.
Then $f$ and $g$ are $\times$-homotopic if and only if they are in
the same path-connected component of $\Hom(G,H)$.  In particular,
the number of $\times$-homotopy classes of maps from $G$ to $H$ is
equal to the number of path components in $\Hom(G,H)$.
\end{prop}

\begin{proof}
Suppose $f,g:G \rightarrow H$ are graph maps such that $f$ and $g$
are in the same component of $\Hom(G,H)$.  Then we can find a path
from $f$ to $g$ in $|\Hom(G,H)|$, which can be approximated as a
finite walk $(f, x_1, x_2, \dots, g)$ on the 1-skeleton.  We claim
that we can extend this to a walk $(f = f_0, x_1, f_1, x_2, f_2,
\dots, f_n = g)$, where each $f_i:G \rightarrow H$ is a graph map
(i.e., $f_i(v)$ consists of a single element for each $v \in V(G)$).

To see this, note that $f \leq x_1$ in $\Hom(G,H)$.  First suppose
that $x_1 \leq x_2$.  Then for each $v \in V(G)$, we choose (by the
choice axiom, say) a single element of $x_1(v)$ to get our map $f_1:
G \rightarrow H$ such that $f_1 \leq x_1 \leq x_2$.  Next suppose
$x_2 \leq x_1$.  If $x_2$ is already a graph map, take $f_1 = x_2$,
and otherwise for each $v \in V(G)$ choose a single element of
$x_2(v)$ to get a map $f_1:G \rightarrow H$.

Now, to get our homotopy, we define a map $F: G \times I_n
\rightarrow H$ by $F(v,i) = f_i(v)$.  Then $F$ is indeed a graph map
since we have an $x_i \in \Hom(G,H)$ such that $f_{i-1}, f_i \leq
x_i$ for each $0 < i \leq n$.  Hence the maps $f=f_0$ and $g=f_n$
are $\times$-homotopic.

For the other direction, suppose that $f,g:G \rightarrow H$ are
distinct maps that are $\times$-homotopic for $n = 1$.  We define a
function $\xi: V(G) \rightarrow 2^{V(H)} \backslash \{\emptyset\}$
by $v \mapsto \{f(v), g(v)\}$.  We claim that $\xi$ is a cell in
$\Hom(G,H)$.  To see this, suppose $v \sim w$ are adjacent vertices of $G$.  Then both
$(f(v),f(w))$ and $(g(v),g(w))$ are edges in $H$ since $f$ and $g$
are graph maps. Also, $(0v, 1w)$ and $(0w,1v)$ are edges in $G
\times I_1$ and since $f$ and $g$ are 1-homotopic this implies that
$(f(v),g(w))$ and $(f(w),g(v))$ are both edges in $H$.  Hence
vertices of $\xi(v)$ are adjacent to vertices of $\xi(w)$ as
desired.  It is clear that both $f$ and $g$ are vertices of $\xi$
and hence we have a path from $f$ to $g$. Now, suppose $f$ and $g$
are $\times$-homotopic for some choice of $n$ and let $F:G \times
I_n \rightarrow H$ be the homotopy.  Let $f_i:G \mapsto H$ be the
graph map given by $v \mapsto F(\iota_i(v))$. Then by induction we
have a path in $\Hom(G,H)$ from $f$ to $f_{n-1}$ and the above
construction gives a path from $f_{n-1}$ to $f_n = g$.
\end{proof}

We end this section with the following observation.

\begin{lemma}
Let $G$ be a graph, $k \leq n$ integers, and let $\iota_k:G
\rightarrow G \times I_n$ denote the graph map given by $\iota(g) =
(g,k)$.  Then for a graph $T$, the induced map $\iota_{k_T}:
\Hom(T,G) \rightarrow \Hom(T, G \times I_n)$ is a homotopy
equivalence.
\end{lemma}

\begin{proof}
Let $i:\Hom(T,G) \times \Hom(T,I_n) \hookrightarrow \Hom(T, G \times
I_n)$ denote the inclusion, a homotopy equivalence by Proposition
\ref{product}.  Let $\phi_k:\Hom(T,G) \rightarrow \Hom(T,G) \times
\Hom(T, I_n)$ denote the inclusion given by $x \mapsto (x,c_k)$,
where $c_k \in \Hom(T,I_n)$ is the constant map sending all elements
of $V(T)$ to $k$.  We note that $\phi_k$ is a homotopy equivalence
by Lemma \ref{contractible}. We then have the following commutative
diagram showing that $\iota_{k_T} = i \circ \phi_k$ is a homotopy
equivalence.

\begin{center}
$\xymatrix{ \Hom(T,G) \ar[r]^{\iota_{k_T}} \ar[d]^{\phi_k}_{\simeq}
& \Hom(T, G \times I_n) \\
\Hom(T,G) \times \Hom(T,I_n) \ar[ur]^{i}_{\simeq} }$
\end{center}
\end{proof}

\section{Homotopy equivalence of graphs}

If $f, g: G \rightarrow H$ are graph maps, the functors obtained by
fixing a graph $T$ in one coordinate of the $\Hom$ complex in each
case provides a pair of topological maps.  For a fixed test graph
$T$, the functor $\Hom(T,?)$ provides the pair of maps $f_T, g_T:
\Hom(T,G) \rightarrow \Hom(T,H)$, while $\Hom(?,T)$ provides the
maps $f^T, g^T:\Hom(H,T) \rightarrow \Hom(G,T)$ (discussed above).  If $f$ and $g$ are
$\times$-homotopic, we can ask how these induced maps are related up
to (topological) homotopy.  It turns out that the induced maps are
homotopic, and in fact provide a characterization of graph
$\times$-homotopy in each case.  More precisely, we have the
following result.

\begin{thm} \label{homotopic}
Let $f,g:G \rightarrow H$ be graph maps.  Then the following are
equivalent:

(1) $f$ and $g$ are $\times$-homotopic.

(2) For every graph $T$, the induced maps $f_T, g_T:\Hom(T,G)
\rightarrow \Hom(T,H)$ are homotopic.

(3) The induced maps $f_G, g_G:\Hom(G,G) \rightarrow \Hom(G,H)$ are
homotopic.

(4) For every graph $T$, the induced maps $f^T, g^T:\Hom(H,T)
\rightarrow \Hom(G,T)$ are homotopic.

(5) The induced maps $f^H, g^H:\Hom(H,H) \rightarrow \Hom(G,H)$ are
homotopic.

\end{thm}

\begin{proof}
We first prove $(1) \Rightarrow (2)$.  Suppose $f,g:G \rightarrow H$
are $\times$-homotopic via a graph map $F:G \times I_n \rightarrow
H$. Then (with notation as above) we have a commutative diagram in
${\mathcal G}$ and, via the functor $\Hom(T,?)$, the induced diagram
in ${\mathcal TOP}$, the category of topological spaces and
continuous maps, of the form:

\medskip

\begin{center}
$\xymatrix{
G \ar[d]_{\iota_0} \ar[dr]^{f} \\
G \times I_n \ar[r]^F & H \\
G \ar [u]^{\iota_n} \ar [ur]_g }$ \hspace {0.5 in} $\xymatrix{
\Hom(T,G) \ar[d]_{\iota_{0_T}} \ar[dr]^{f_T} \\
\Hom(T, G \times I_n) \ar[r]^{F_T} & \Hom(T, H) \\
\Hom(T,G) \ar[u]^{\iota_{n_T}} \ar [ur]_{g_T} }$
\end{center}

\medskip

Now, $\Hom(T,I_n)$ is path connected (contractible) by Lemma
\ref{contractible}. Let $\gamma:I=[0,1]\rightarrow \Hom(T,I_n)$ be a
path such that $\gamma(0) = c_0$ and $\gamma(1) = c_n$ (where again
$c_i \in \Hom(T,I_n)$ is the constant map sending all vertices of
$T$ to $i$).  Let $j_i:\Hom(T,G) \rightarrow \Hom(T,G) \times I$ be
the (topological) map given by $(\id,i)$.  We then obtain the
following diagram in ${\mathcal TOP}$ (where $(T,G) = \Hom(T,G)$,
etc.):

\medskip

\begin{center}
$\xymatrix{
& (T,G) \ar[dl]_{j_0} \ar[dr]^{\iota_{0_T}} \\
(T,G) \times I \ar[r]^{id \times \gamma} & (T,G) \times (T, I_n) \ar @{^{(}->}[r] & (T, G \times I_n) \ar[r]^{F_T} & (T, H) \\
& (T,G) \ar[ul]^{j_1} \ar[ur]_{\iota_{n_T}} }$
\end{center}

\medskip

We claim that this diagram commutes.  To see this, suppose $\alpha
\in \Hom(T,G)$. Then for all $t \in V(T)$ we have
$\iota_{0_T}(\alpha)(t) = \{\iota_0(x):x \in \alpha(t)\} = \{(x,0):x
\in \alpha(t)\} \in \Hom(T,G) \times \Hom(T,I^n)$, so that
$\iota_{0_T}(\alpha) = (\alpha, c_0)$.  On the other hand, $(id
\times \gamma)(j_0)(\alpha) = (id \times \gamma)(\alpha, 0) =
(\alpha, c_0)$.  The bottom square is similar.

Now, let $\Phi:\Hom(T,G) \times I \rightarrow \Hom(T,H)$ be the
composition from above.  We have that $\Phi \circ j_0 = F_T \circ
\iota_{0_T} = f_T: \Hom(T,G) \rightarrow \Hom(T,H)$ and similarly
$\Phi \circ j_1 = g_T$, so that $f_T$ and $g_T$ are homotopic.

The implication $(2) \Rightarrow (3)$ is clear, so we next turn to
$(3) \Rightarrow (1)$.  For this, suppose $f,g:G
\rightarrow H$ are not $\times$-homotopic.  Then $f$ and $g$ are in
different path components of $\Hom(G,H)$ by Proposition
\ref{component}.  We claim that the induced maps $f_G, g_G:
\Hom(G,G) \rightarrow \Hom(G,H)$ are also not homotopic.  To obtain
a contradiction, suppose they are and let $\Phi:\Hom(G,G) \times I
\rightarrow \Hom(G,H)$ be a (topological) homotopy between them.
Note that if $\id \in \Hom(G,G)$ is the identity map, then $f_G(\id)
= f$ and $g_G(\id) = g$ since, for instance, we have $f_G(\id)(x) =
\{f(y):y \in \id(x)\} = \{f(y):y \in \{x\}\} = \{f(x)\}$ for all $x
\in V(G)$.  So then the restriction $\Phi|_{\{\id\} \times
I}:\Hom(G,G) \times I \rightarrow \Hom(G,H)$ gives a path in
$\Hom(G,H)$ from $f$ to $g$, a contradiction.

We next prove $(1) \Rightarrow (4)$.  Again, suppose $f,g:G
\rightarrow H$ are $\times$-homotopic via $F:G \times I_n
\rightarrow H$. Then this time we have the commutative diagram in
${\mathcal G}$ and the induced diagram in ${\mathcal TOP}$ of the
form:

\medskip

\begin{center}
$\xymatrix{
G \ar[d]_{\iota_0} \ar[dr]^{f} \\
G \times I_n \ar[r]^F & H \\
G \ar[u]^{\iota_n} \ar [ur]_g }$ \hspace {0.5 in} $\xymatrix{
& \Hom(G,T) \\
\Hom(H,T) \ar[r]^{F^T} \ar[ur]^{f^T} \ar [dr]_{g^T} & \Hom(G \times I_n, T) \ar[u]_{\iota_{0}^T} \ar[d]^{\iota_{n}^T}\\
& \Hom(G,T) }$
\end{center}

\medskip

To show that $f^T$ and $g^T$ are homotopic, we will find a map
$\Psi: \Hom(H,T) \rightarrow \Hom(G,T)^I$ such that $p_0 \Psi = f^T$
and $p_1 \Psi = g^T$.  First, we define a map $\varphi:\Hom(I_n,
T^G) \times \{0, \frac{1}{n}, \frac{2}{n}, \dots, 1\} \rightarrow
\Hom({\bf 1},T^G)$ via $\varphi(\alpha,\frac{i}{n})(v) = \alpha(i)$
for $v \in {\bf 1}$, $\alpha \in \Hom(I_n, T^G)$, and $0 \leq i \leq
n$.  This extends to a map $\varphi:\Hom(I_n,T^G) \times I
\rightarrow \Hom({\bf 1}, T^G)$ since the maps
$\varphi_j:\Hom(I_n,T^G) \rightarrow \Hom({\bf 1}, T^G)$ are all
homotopic for $0 \leq j \leq n$ (recall $\iota_j:{\bf 1} \rightarrow
I_n$ induces a homotopy equivalence).  Let $\tilde
\varphi:\Hom(I_n,T^G) \rightarrow \Hom({\bf 1}, T^G)^I$ be the
adjoint map.   Next, from the above proposition, we have a map
$\psi:\Hom({\bf 1}, T^G) \rightarrow \Hom(G,T)$ that is a homotopy
inverse to the inclusion $\Hom({\bf 1} \times G, T) \rightarrow
\Hom({\bf 1}, T^G)$.  Let $\tilde \psi: \Hom({\bf 1}, T^G)^I
\rightarrow \Hom(G,T)^I$ be the induced map on the path spaces.
Define $\Phi:\Hom(I_n, T^G) \rightarrow \Hom(G,T)^I$ by the
composition $\Phi = \tilde \varphi \tilde \psi$. Finally, we get the
desired map $\Psi$ as the horizontal composition in the commutative
diagram below.

\begin{center}
$\xymatrix{
& & (G,T) \\
(H,T) \ar[r]^{F^T} & (G \times I_n, T) \ar [ur]^{\iota_{0_T}}
\ar[dr]_{\iota_{n_T}} \ar @{^{(}->}[r] & (I_n, T^G) \ar[r]^{\Phi} &
(G,T)^I \ar[ul]_{p_0} \ar[dl]^{p_1} \\
& & (G,T) }$

\end{center}

The implication $(4) \Rightarrow (5)$ is again clear, and so we are
left with only $(5) \Rightarrow (1)$.  For this, suppose $f,g:G
\rightarrow H$ are not $\times$-homotopic, so that $f$ and $g$ are
in different path components of $\Hom(G,H)$.  We claim that the
induced maps $f^H, g^H: \Hom(H,H) \rightarrow \Hom(G,H)$ are not
homotopic. Suppose not, so that we have $f^H, g^H:\Hom(H,H)
\rightarrow \Hom(G,H)$ are homotopic via a (topological) map
$\Phi:\Hom(H,H) \times I \rightarrow \Hom(G,H)$.  Here note that if
$\id \in \Hom(H,H)$ is the identity map, then $f^H(\id) = f$ and
$g^H(\id) = g$ since, for instance, $f^H(\id)(x) = \id(f(x)) =
f(x)$.  Hence the restriction $\Phi|_{\{\id\} \times I}:\Hom(H,H)
\times I \rightarrow \Hom(G,H)$ gives a path in $\Hom(G,H)$ from $f$
to $g$, a contradiction.  The result follows.
\end{proof}

The notion of $\times$-homotopy of graph maps provides a natural
candidate for the notion of $\times$-homotopy equivalence of graphs.  Again,
this has several equivalent formulations, which we establish next.

\begin{thm} \label{equivalent}
Let $f: G \rightarrow H$ be maps of graphs.  Then the following are
equivalent.

(1)  There exists a map $g:H \rightarrow G$ such that $f \circ g
\simeq_\times \id_H$ and $g \circ f \simeq_\times \id_G$ (call $g$ a
\textit{homotopy inverse} to $f$).

(2)  For every graph $T$, the induced map $f_T: \Hom(T,G)
\rightarrow \Hom(T,H)$ is a homotopy equivalence.

(3)  For every graph $T$, the induced map $(f_{T})_0 :
\pi_0 \big(\Hom(T,G)\big) \rightarrow \pi_0 \big(\Hom(T,H)\big)$ is an isomorphism (bijection).

(4)  For every graph $T$, the induced map $f_T: [T,G]_{\times}
\rightarrow [T,H]_{\times}$ is a bijection.

(5)  The maps $f_G:\Hom(G,G) \rightarrow \Hom(G,H)$ and
$f_H:\Hom(H,G) \rightarrow \Hom(H,H)$ both induce isomorphisms on
the path components.

(6)  For every graph $T$, the induced map $f^T: \Hom(H,T)
\rightarrow \Hom(G,T)$ is a homotopy equivalence.

(7)  The maps $f^G:\Hom(H,G) \rightarrow \Hom(G,G)$ and
$f^H:\Hom(H,H) \rightarrow \Hom(G,H)$ both induce isomorphisms on
path components.

\begin{center}
$\xymatrix{
& \Hom(G,H) \\
\Hom(G,G) \ar[ur]^{f_G} && \Hom(H,H) \ar[ul]^{f^H} \\
& \Hom(H,G) \ar[ul]_{f^G} \ar[ur]_{f_H} }$

\end{center}
\end{thm}

\begin{proof}
For $(1) \Rightarrow (2)$, $g_T$ is a homotopy inverse by Theorem
\ref{homotopic}.

$(2) \Rightarrow (3)$ is clear.

$(3) \iff (4)$ follows from Proposition \ref{component}.

$(3) \Rightarrow (5)$ is clear.

For $(5) \Rightarrow (1)$, we assume $(f_{H})_0 : \pi_0 \big(\Hom(H,G)\big)
\rightarrow \pi_0 \big(\Hom(H,H)\big)$ is an isomorphism.  Let $\phi$ be its
inverse and let $(\id_H)_0$ denote the connected component of
$\id_H$ in $\Hom(H,H)$.  Let $g \in \phi \big((\id_H)_0 \big)$ be a vertex
of $\Hom(H,G)$ (i.e., a graph map).  We claim that $g$ satisfies the
conditions of (1). To see this note that $\big((f_H)_0 \phi \big)\big((\id_H)_0 \big)
= (\id_H)_0$ and since $g \in \phi \big((\id_H)_0\big)$ we have that $fg =
f_H(g)$ is in the same component as $\id_H$ in $\Hom(H,H)$.  Hence
$fg \simeq_\times \id_H$, as desired.  A similar consideration of
the isomorphism $(f_{G})_0 : \pi_0 \big(\Hom(G,G)\big) \rightarrow
\pi_0 \big(\Hom(G,H) \big)$ shows that $gf \simeq_\times \id_G$.

For $(1) \Rightarrow (6)$, $g^T$ again provides the inverse by
Theorem \ref{homotopic}.

$(6) \Rightarrow (7)$ is clear.

Finally, we check $(7) \Rightarrow (1)$.  For this we assume
$(f^G)_0: \pi_0 \big(\Hom(H,G)\big) \rightarrow \pi_0 \big(\Hom(G,G)\big)$ is an
isomorphism.  Let $\psi$ be the inverse and let $(\id_G)_0$ denote
the connected component of $\id_G$.  Let $g \in \psi \big((\id_G)_0 \big)$ be
a graph map $g:H \rightarrow G$.  We claim that $g$ satisfies the
conditions that we need.   Note that $\big((f^G)_0 \psi \big) \big((\id_G)_0 \big) =
(\id_G)_0$ and $f^G(g) = gf$, and hence $gf \simeq_\times \id_G$.
Similarly we get $fg \simeq_\times \id_H$ and the result follows.
\end{proof}

\begin{defn}
A graph map $f:G \rightarrow H$ is called a \textit{$\times$-homotopy
equivalence} (or simply homotopy equivalence) if it satisfies any of the above conditions.  Homotopy
equivalence of graphs is an equivalence relation, and we let $[G]$
denote the homotopy equivalence class of $G$.
\end{defn}

Aside from certain qualitative similarities, homotopy equivalences
of graphs satisfy many of the formal properties enjoyed by
equivalences in any abstract homotopy theory, \cite{HOV99} and \cite{QUI67}.  We close this section with a couple of observations along these lines.

\begin{defn}
Let ${\mathcal M}$ be a class of maps in a category ${\mathcal C}$.
${\mathcal M}$ is said to satisfy the \textit{2 out of 3 property}
if, for all maps $f$ and $g$, whenever any two of $f,g,gf$ are in
${\mathcal M}$, then so is the third.
\end{defn}

\begin{lemma}
Homotopy equivalences of graphs satisfy the 2 out of 3 property.
\end{lemma}

\begin{proof}
Let $f:X \rightarrow Y$ and $g:Y \rightarrow Z$ be maps of graphs,
and let $T$ be a graph.  We will be considering the following
diagrams.

\begin{center}

$\xymatrix{ \Hom(T,X) \ar@<1ex>[r]^{f_T} & \Hom(T,Y)
\ar@<1ex>[l]^{a_T} }$

$\xymatrix{ \Hom(T,Y) \ar@<1ex>[r]^{g_T} & \Hom(T,Z)
\ar@<1ex>[l]^{b_T} }$

$\xymatrix{ \Hom(T,X) \ar@<1ex>[r]^{gf_T} & \Hom(T,Y)
\ar@<1ex>[l]^{c_T} }$.

\end{center}

First suppose $f$ and $g$ are both homotopy equivalences, with
homotopy inverse maps $a:Y \rightarrow X$ and $b:Z \rightarrow Y$
respectively.  We claim $ab$ is the homotopy inverse to $fg$.  To
see this, note that $(abgf)_T = a_T b_T g_T f_T \simeq a_T f_T
\simeq \id_X$.  Similarly, $(gfab)_T \simeq \id_Z$, so that $gf$ is
a homotopy equivalence.

Next suppose that $f$ and $gf$ are homotopy equivalences, and let
$c: Z \rightarrow X$ be the homotopy inverse to $gf$.  We claim
$fc:Z \rightarrow Y$ is the homotopy inverse to $g$.  For this we
compute $(gfc)_T = g_T f_T c_T \simeq \id_Z$ and $(fcg)_T = f_T c_T
g_T \simeq f_T c_T g_T f_T a_T \simeq f_T a_T \simeq \id_Y$.  We
conclude that $g$ is a homotopy equivalence.

Finally, we claim that if $g$ and $gf$ are homotopy equivalences
then $cg:Y \rightarrow X$ is the homotopy inverse to $f$.  This
follows from the fact that $(fcg)_T = f_T c_T g_T \simeq b_T g_T f_T
c_T g_T \simeq b_T g_T \simeq \id_Y$ and also $(cgf)_T = c_T g_T f_T
\simeq \id_Z$.
\end{proof}

\begin{defn}
Let $g:G \rightarrow H$ be a map in a category ${\mathcal C}$.
Recall that $f$ is a \textit{retract} of $g$ if there is a
commutative diagram of the following form,

\begin{center}
$\xymatrix{
X \ar[r]^\alpha \ar[d]_f & G \ar[r]^\gamma \ar[d]_g & X \ar[d]^f  \\
Y \ar[r]_\beta & H \ar[r]_\delta & Y }$
\end{center}

\noindent where the horizontal composites are identities.

\end{defn}

\begin{lemma} Homotopy equivalences of graphs are closed under retracts.
\end{lemma}

\begin{proof}
Suppose $g$ is a homotopy equivalence.  Then for every graph $T$ we
have the diagram,

\begin{center}
$\xymatrix{
\Hom(T,X) \ar[r]^{\alpha_T} \ar[d]_{f_T} & \Hom(T,G) \ar[r]^{\gamma_T} \ar[d]_{g_T}^{\simeq} & \Hom(T,X) \ar[d]^{f_T}  \\
\Hom(T,Y) \ar[r]_{\beta_T} & \Hom(T,H) \ar[r]_{\delta_T} & \Hom(T,Y)
}$
\end{center}

\noindent with $g_T:\Hom(T,G) \rightarrow \Hom(T,H)$ a homotopy
equivalence. We consider the induced maps on homotopy groups.  Since
$\gamma_T \alpha_T = \id$, we have that $(\alpha_T)_*$ is injective
and hence so is $(f_T)_*$, since $(\beta_T)_* (f_T)_* = (g_T)_*
(\alpha_T)_*$ is injective.  Similarly, since $\delta_T \beta_T
=\id$, we have that $(\delta_T)_*$ is surjective and hence so is
$(f_T)_*$. We conclude that $f_T$ induces an isomorphism on all
homotopy groups and hence $f_T$ is a homotopy equivalence on the
$CW$-type $\Hom$ complexes.
\end{proof}

\section{Foldings, stiff graphs, and dismantlable graphs}

In this section we investigate some further properties and consequences of $\times$-homotopy of graphs.  The relevant operation in this context will that of a graph \textit{folding}, which we will see is closely related to $\times$-homotopy.

\begin{defn}

Let $u$ and $v$ be vertices of a graph $G$ satisfying $N(v)
\subseteq N(u)$.  Then the map $f:G \rightarrow G \backslash v$
given by $f(x) = x$,  $x \neq v$, and $f(v) = u$, is called a
\textit{folding} of $G$ at the vertex $v$. Similarly, the inclusion
$i:G \backslash v \rightarrow G$ is called an \textit{unfolding}
(see Figure 10).
\end{defn}

\begin{center}
\epsfig{file=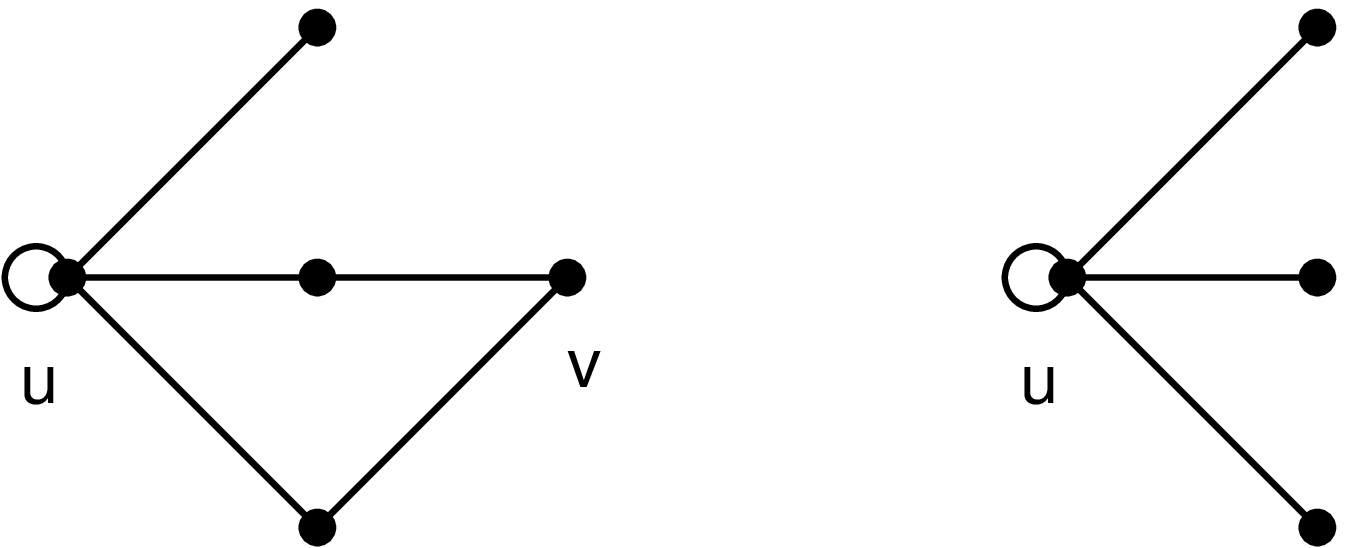, height=1.0 in, width = 2.5 in}

{Figure 10: The graph $G$ and the folded graph $G \backslash v$}

\end{center}

In the original papers regarding $\Hom$ complexes (see for example
\cite{BKcom}), it was shown that foldings in the first coordinate of
the $\Hom$ complex preserved homotopy type.  For some time it was an
open question whether the same was true in the second coordinate of
the $\Hom$ complex.  Kozlov investigated this question in the papers
\cite{Kfold} and \cite{Kcol}, and showed that indeed this was the
case.

\begin{prop}[Kozlov] \label{folds}
If $G$ and $H$ are graphs, and $u$ and $v$ are vertices of $G$ such
that $N(v) \subseteq N(u)$, then the folding and unfolding maps
induce inclusions of strong deformation retracts

\begin{center}

$\xymatrix{ \Hom(G \backslash v,H) \ar@{^{(}->}[r]_{\simeq}^{f^H} &
\Hom(G,H),}$
\hspace{.2 in}
$\xymatrix{ \Hom(H,G \backslash v) \ar@{^{(}->}[r]_{\simeq}^{i_H} &
\Hom(H,G). }$

\end{center}
\end{prop}

In fact, Kozlov exhibits these deformation retracts as closure maps
on the levels of the posets, which he shows preserve the
\textit{simple} homotopy type of the associated simplicial complex (we refer to \cite{Kozsimp} for necessary definitions).
 We note that although Kozlov deals only with the situation of finite
$H$, his proof extends to the case of arbitrary $H$.  In Sections 5
and 6 of this paper we see the further importance of folds in the
context of the $\Hom$ complex.

\begin{rem}
We can apply Theorem \ref{equivalent} to obtain the following alternate proof of
one part of Proposition \ref{folds}.  As we
mentioned, it was previously known that if $G \rightarrow H = G \backslash
\{v\}$ is a folding, then $f^T:\Hom(H,T) \rightarrow \Hom(G,T)$ is a
homotopy equivalence for all $T$.  We can then apply $(6)
\Rightarrow (2)$ in Theorem \ref{equivalent} to conclude that
$f_T:\Hom(T,G) \rightarrow \Hom(T,H)$ is also a homotopy
equivalence, and hence `folds in the second coordinate' also
preserve homotopy type of $\Hom$ complexes.  Our theorem also provides some insight into the symmetry involved in the two entries of the $\Hom$ complex.
\end{rem}

\subsection{Stiff graphs}

If $f:G \rightarrow \tilde G$ is a map realized by a sequence of foldings and unfoldings, then
$f_T:\Hom(T,G) \rightarrow \Hom(T,\tilde G)$ is a homotopy
equivalence for all $T$, and hence $G$ and $\tilde G$ are homotopy
equivalent.  One can then consider the case when $G$ has no more
foldings available.  From \cite{HN04} we have the following notion.

\begin{defn}
A graph $G$ is called \textit{stiff} if there does not exist a
pair of distinct vertices $u, v \in V(G)$ such that $N(v) \subseteq
N(u)$.
\end{defn}

\begin{lemma}
Suppose $G$ is a stiff graph.  Then the identity map $\id_G$ is an
isolated point in the realization of $\Hom(G,G)$.
\end{lemma}

\begin{proof}
If not, then we have some $\alpha \in \Hom(G,G)$ such that $x \in
\alpha(x)$ for all $x \in V(G)$, and such that $\{v,w\} \subseteq
\alpha(v)$ for some $v \neq w$.  Since $G$ is stiff we have some
vertex $x \in V(G)$ such that $x \in N(v) \backslash N(w)$.  But
then since $x \in \alpha(x)$ we need $x$ to be adjacent to $w$ (to
satisfy the conditions of $\Hom$), a contradiction.
\end{proof}

\begin{prop} \label{stiff}
If $G$ and $H$ are both stiff graphs, then $G$ and $H$ are homotopy
equivalent if and only if they are isomorphic.
\end{prop}

\begin{proof}
Sufficiency is clear.  For the other direction, suppose $f:G
\rightarrow H$ is a homotopy equivalence with inverse $g:H
\rightarrow G$.  Then $gf$ is $\times$-homotopic to the identity
$\id_G$, so that $gf$ and $\id_G$ are in the same component of
$\Hom(G,G)$ by Proposition \ref{component}.  But then $gf = \id_G$
since $G$ is stiff. Similarly we get $fg = \id_H$, so that $f$ is an
isomorphism.
\end{proof}

From this it follows that if $G$ and $H$ are finite graphs and $f:G
\rightarrow H$ is a homotopy equivalence, then one can fold both
graphs to their unique (up to isomorphism) stiff subgraphs $\tilde
G$ and $\tilde H$ and get an isomorphism $\tilde G = \tilde H$.
However, in general one cannot make these foldings commute with
the map $f$, as the next example illustrates.

\begin{example}
Let $G$ be the graph with 5 vertices $V(G) = \{1,2,3,4,5\}$ and
edges $E(G) = \{11,12,15,22,23,33,35,34,44,45\}$ (see Figure 11).
Let $f:{\bf 1} \rightarrow G$ be the map that maps ${\bf 1} \mapsto
4$.

\begin{center}
\epsfig{file=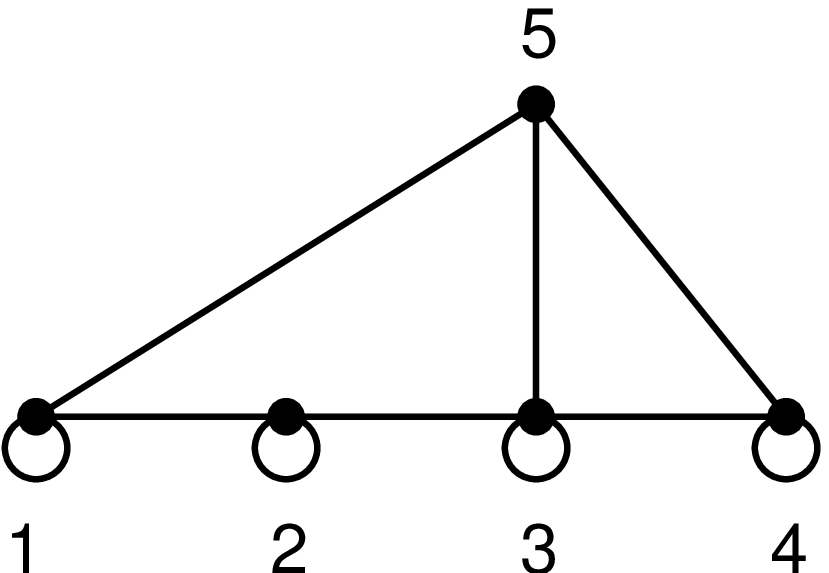, height=1.0 in, width = 1.4 in}

{Figure 11: The graph $G$}
\end{center}

\noindent We note that $G$ is foldable to a looped vertex ${\bf 1}$ (and hence homotopy equivalent to $G$),
but cannot be folded to $\im(f)$ by a sequence of foldings and unfoldings.
\end{example}

\begin{question}
Suppose $G$ and $H$ are (finite) graphs and $f:G \rightarrow H$ is a
homotopy equivalence.  Under what circumstances can $f$ be factored as a sequence of foldings
and unfoldings?
\end{question}

\noindent Note that an answer to this question would
yield another characterization of homotopy equivalence to the list
in Theorem \ref{equivalent}, under the relevant conditions on $G$
and $H$.

\textit{(8) The graph map $f:G \rightarrow H$ can be factored as a
sequence of foldings and unfoldings.}

\subsection{Dismantlable graphs}

As in \cite{HN04}, a finite graph $G$ is called \textit{dismantlable} if it can be
folded down to ${\bf 1}$.  Note that $G$ is dismantlable if
\textit{any} sequence of foldings of $G$ down to its stiff subgraph
results in the looped vertex ${\bf 1}$.  Dismantlable graphs have
gained some attention in the recent papers of Brightwell and Winkler
(see \cite{BW00} and \cite{BW04}), where they are related to the
uniqueness of Gibbs measure on the set of homomorphisms between two
graphs.  We can apply the results of Theorem \ref{equivalent} to
obtain the following characterizations of dismantlable graphs.

\begin{prop}
Suppose $G$ is a finite graph, and let $f:G \rightarrow {\bf 1}$ be
the unique map.  Then the following are equivalent:

(0) $G$ is dismantlable.

(1) There exists a map $g:{\bf 1} \rightarrow G$ such that $fg
\simeq_\times \id_{\bf 1}$ and $gf \simeq_\times \id_G$.

(2) For every graph $T$, the map $f_T:\Hom(T,G) \rightarrow
\Hom(T,{\bf 1})$ is a homotopy equivalence.

(2a) For every graph $T$, $\Hom(T,G)$ is contractible.

(3) For every graph $T$, $\Hom(T,G)$ is connected.

(4) For every graph $T$, the set $[T,G]_\times$ consists of a single
homotopy class.

(5) $G$ has at least one looped vertex and $\Hom(G,G)$ is connected.

(6) The map $f^G:\Hom({\bf 1},G) \rightarrow \Hom(G,G)$ induces an
isomorphism on path components.
\end{prop}

\begin{proof}

$(1) \Rightarrow (2)$ is a special case of Theorem \ref{equivalent}
(with $H = {\bf 1}$), and $(2) \Rightarrow (2a)$ since $\Hom(T,{\bf
1})$ is contractible for all $T$.

$(2a) \Rightarrow (3)$ is clear, and the sequence of equivalences
$(3) \iff (4) \iff (1)$ is again a special case of Theorem
\ref{equivalent}.

The implication $(3) \Rightarrow (5)$ is clear.  For $(5)
\Rightarrow (1)$, we assume that $v \in V(G)$ is a looped vertex,
and that $\Hom(G,G)$ is connected.  Let $g:{\bf 1} \rightarrow G$ be
the graph map given by ${\bf 1} \rightarrow v$. We claim that $g$
satisfies the conditions of (1). First, we have $fg = \id_{\bf 1}$.
Also, since $\Hom(G,G)$ is path connected, we have that $gf:G
\rightarrow G$ is in the same path component as the identity
$\id_G$. Hence $gf \simeq_\times \id_G$, and so $g$ is the desired
graph map.

Finally, $(6) \iff (1)$ is another special case of Theorem
\ref{equivalent}. Here note that $f^{\bf 1}:\Hom({\bf 1}, {\bf 1})
\rightarrow \Hom(G, {\bf 1})$ is always an isomorphism.

It only remains to show $(0) \iff (3)$.  If $G$ is foldable to a
looped vertex then Proposition \ref{folds} implies that $\Hom(T,G) \simeq \Hom(T, {\bf
1})$; the latter space is a point (and hence connected) for all $T$.
For the other direction, we suppose $\Hom(T,G)$ is connected for all
graphs $T$. The unique map $G \rightarrow {\bf 1}$ gives a bijection
$\pi_0 \big((\Hom(T,G)\big) \rightarrow \pi_0 \big(\Hom(T,{\bf 1}) \big)$ for all $T$,
and hence $G$ and ${\bf 1}$ are homotopy equivalent.  So then if $G$
is stiff, we have that $G$ is isomorphic to ${\bf 1}$ by Proposition
\ref{stiff}. Otherwise we perform folds to reduce the number of
vertices and use induction on $|V(G)|$.
\end{proof}

\section{Other internal homs and $A$-theory}

In this last section we investigate other notions of graph homotopy
that arise under considerations of different internal hom
structures. One such homotopy theory (associated to the cartesian
product) recovers the \textit{$A$-theory} of graphs as defined in \cite{BL}.

Recall that in our construction of $\times$-homotopy, we relied on
the fact that the categorical product has the looped vertex at its
unit, and also possesses an internal hom (exponential) construction.
This meant that graph maps from $G$ to $H$ were encoded by the
looped vertices in the graph $H^G$, and two maps $f,g:G \rightarrow
H$ were considered $\times$-homotopic if one could walk from $f$ to
$g$ along a path composed of other graph maps.

Hence, in the general set-up we will be interested in monoidal category structures on the category of graphs that have the
looped vertex as the unit element (this just means that we have an associative bifunctor $\otimes:{\mathcal G} \times {\mathcal G} \rightarrow {\mathcal G}$), together with an internal hom for
that structure. Recall that \textit{having an internal hom} means
that the set valued functor $T \mapsto {\mathcal G}(T \otimes G, H)$
is representable by an object of ${\mathcal G}$, which we will
denote by $H^G$.  We then have $T \mapsto {\mathcal G}(T \otimes G,
H) = {\mathcal G}(T,H^G)$. Since we require the looped vertex (which
we denote by ${\bf 1}$) to be the unit we also get ${\mathcal
G}(G,H) = {\mathcal G}({\bf 1} \otimes G, H) = {\mathcal G}({\bf 1},
H^G)$, so that $H^G$ is a graph with the looped vertices as
precisely the set of graph maps $G \rightarrow H$.  A pair of graph
maps $f$ and $g$ will then be considered homotopic in this context
if, once again, we can find a (finite) path from $f$ to $g$ along
looped vertices.

One such product of interest is the \textit{cartesian product}; we
recall its definition below.

\begin{defn}
For graphs $A$ and $B$, the \textit{cartesian product} $A \square B$
is the graph with vertex set $V(A) \times V(B)$ and adjacency given
by $(a,b) \sim (a^\prime, b^\prime)$ if either $a \sim a^\prime$ and
$b = b^\prime$, or $a = a^\prime$ and $b \sim b^\prime$ (see Figure
12).
\end{defn}

\begin{center}
\epsfig{file=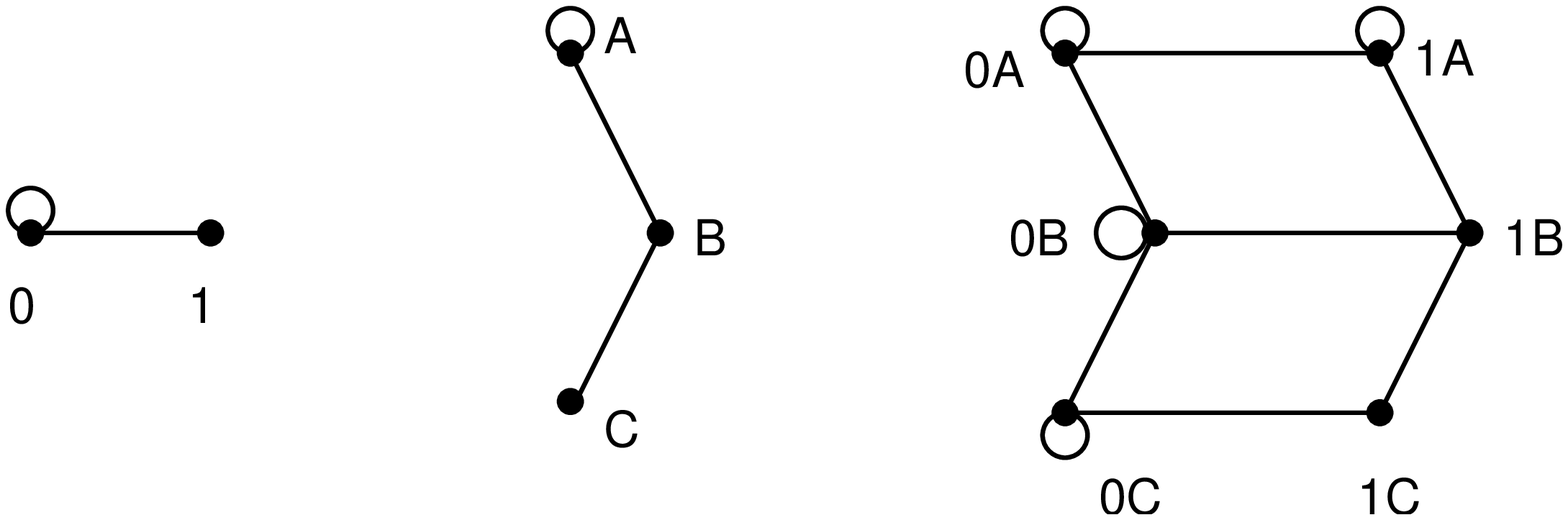, height=1.2 in, width = 4 in}

{Figure 12: The graphs $A$, $B$, and $A \square B$}
\end{center}

One can check that the cartesian product gives the category of
graphs the structure of a monoidal category with a (unlooped) vertex
as the unit element. We next claim that the cartesian product also
has an internal hom; we first define the functor that will serve as
its right adjoint.

\begin{defn}
For graphs $A$ and $B$, the \textit{cartesian exponential graph}
$B^A$ is the graph with vertex set $\{f:A \rightarrow B\}$ the set
of all \textit{graph} maps, with adjacency given by $f \sim
f^\prime$ if $f(a) \sim f^\prime(a)$ for all $a \in A$ (see Figure
13).
\end{defn}

\begin{center}
\epsfig{file=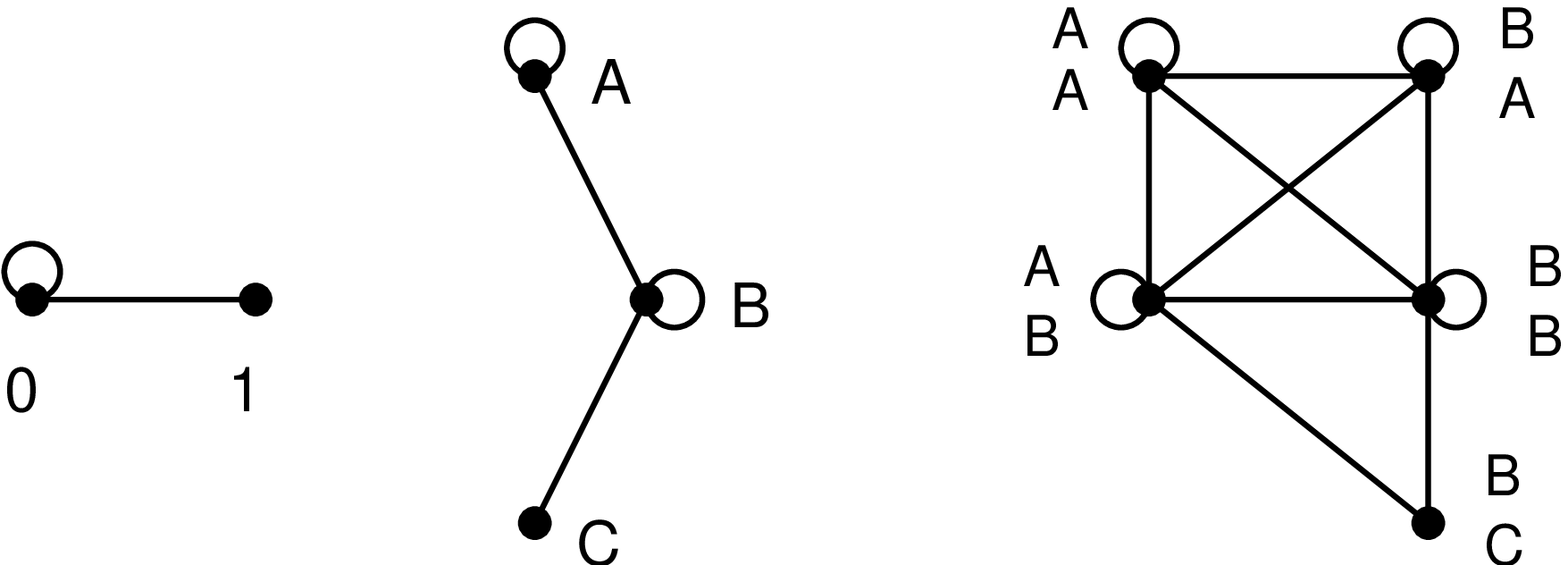, height=1.1 in, width = 3.5 in}

{Figure 13: The graphs $A$, $B$, and $B^A$}
\end{center}

Our next result shows that this exponential construction indeed
provides the right adjoint for the cartesian product defined above.

\begin{lemma} \label{cartadjoint}
For graphs $A,B,C$, there is a natural bijection $\Phi:{\mathcal
G}(A \square B, C) \rightarrow {\mathcal G}(A, C^B)$ given by the
cartesian exponential graph.
\end{lemma}

\begin{proof}
Given $f \in {\mathcal G}(A \square B, C)$, and $a \in V(A)$, $b \in
V(B)$, we define $\Phi(f)(a)(b) = f(a,b)$.  We first verify that
$\Phi(f)(a)$ is a graph map, so that $\Phi(f)(a) \in C^B$.  For
this, suppose $b \sim b^\prime$ are adjacent vertices of $B$.  Then
we have $(a, b) \sim (a, b^\prime)$ in $A \square B$ and hence
$f(a,b) \sim f(a,b^\prime)$ as desired.

Next we verify that $\Phi(f)$ is a graph map.  For this suppose $a
\sim a^\prime$ are adjacent vertices of $A$. Then, once again,
$(a,b) \sim (a^\prime,b)$ in $A \square B$ for all $b \in V(B)$.
Hence $\Phi(f)(a)(b) = f(a,b)$ is adjacent to $\Phi(f)(a^\prime)(b)
= f(a^\prime,b)$ for all $b \in V(B)$, so that $\Phi(f)(a) \sim
\Phi(f)(a^\prime)$.

To see that $\Phi$ is a bijection, we construct an inverse
$\Psi:{\mathcal G}(A, C^B) \rightarrow {\mathcal G}(A \times B, C)$
via $\Psi(g)(a,b) = g(a)(b)(g)$ for all $g \in {\mathcal G}(A,
C^B)$.  One checks that $\Psi$ is well defined and an inverse to
$\Phi$.
\end{proof}

Recall that a \textit{reflexive} graph is a graph with loops on each
vertex, and that a map between reflexive graphs is just a map of the
underlying graphs.  The cartesian product of two reflexive graphs is
once again reflexive, and hence the cartesian product gives the
category ${\mathcal G}^\circ$ of reflexive graphs the structure of a
monoidal category with the looped vertex {\bf 1} as the unit
element.

Also, if $A$ and $B$ are both reflexive, then all vertices of $B^A$
are looped (so that $B^A$ is indeed a reflexive graph).  Hence we
have a graph $B^A$ whose looped vertices are precisely the graph
maps $B \rightarrow A$.  The map $\Phi$ described above then gives a
bijection ${\mathcal G}^\circ (A \square B, C) \simeq {\mathcal
G}^\circ (A, C^B)$.

In some recent papers (see for example \cite{BBLL} and \cite{BL}), a
homotopy theory called $A$-theory has been developed as a way to
capture `combinatorial holes' in simplicial complexes. The
definition can be reduced to a construction in graph theory, applied
to a certain graph associated to the simplicial complex in question.
It turns out that $A$-theory of graphs fits nicely into the set-up
that we have described, where the homotopy theory is associated to
the cartesian product in the category of reflexive graphs.  We
recall the definition of $A$-homotopy of graph maps and $A$-homotopy
equivalence of graphs (as in \cite{BBLL}).

\begin{defn}
Let $f,g:(G,x) \rightarrow (H,y)$ be a pair of based maps of
reflexive graphs. Then $f$ and $g$ are said to be
\textit{$A$-homotopic}, denoted $f \simeq_A g$, if there is an
integer $n \geq 1$ and a graph map $\varphi: G \square I_n
\rightarrow H$ such that $\varphi(?,0) =f$ and $\varphi(?,n) = g$,
and such that $\varphi(x,i) = y$ for all $i$.

We call $(G,x)$ and $(H,y)$ $A$-homotopy equivalent if there exist
based maps $f:G \rightarrow H$ and $g:H \rightarrow G$ such that $gf
\simeq_A \id_G$ and $fg \simeq_A \id_H$.

\end{defn}

Using the adjunction of Lemma \ref{cartadjoint}, we see that an
$A$-homotopy between two based maps of reflexive graphs $f,g:G
\rightarrow H$ is the same thing as a map $\tilde \varphi:I_n
\rightarrow H^G$ with $\tilde \varphi(0) = f$ and $\tilde \varphi(n)
= g$, or in other words a path from $f$ to $g$ along looped vertices
in the based version of the (cartesian) exponential graph $H^G$.
This places the $A$-theory of graphs into the general set-up
described above.  In \cite{BBLL} the authors seek a topological space whose (ordinary) homotopy groups recover the $A$-theory groups of a given graph, and the analogous question in the context of $\times$-homotopy is investigated in \cite{DocGro}.

\bibliographystyle{plain}
\bibliography{litgraph}

\begin{thebibliography}{10}

\bibitem{BBLL}
Eric Babson, H{\'e}l{\`e}ne Barcelo, Mark de~Longueville, and Reinhard
  Laubenbacher.
\newblock Homotopy theory of graphs.
\newblock {\em J. Algebraic Combin.}, 24(1):31--44, 2006.

\bibitem{BKcom}
Eric Babson and Dmitry~N. Kozlov.
\newblock Complexes of graph homomorphisms.
\newblock {\em Israel J. Math.}, 152:285--312, 2006.

\bibitem{BKpro}
Eric Babson and Dmitry~N. Kozlov.
\newblock {Proof of the Lov\'{a}sz Conjecture}.
\newblock {\em Annals of Mathematics}, 165(3):965--1007, 2007.

\bibitem{BL}
H{\'e}l{\`e}ne Barcelo and Reinhard Laubenbacher.
\newblock Perspectives on {$A$}-homotopy theory and its applications.
\newblock {\em Discrete Math.}, 298(1-3):39--61, 2005.

\bibitem{Bjo95}
A.~Bj{\"o}rner.
\newblock Topological methods.
\newblock In {\em Handbook of combinatorics, Vol.\ 1,\ 2}, pages 1819--1872.
  Elsevier, Amsterdam, 1995.

\bibitem{BW00}
Graham~R. Brightwell and Peter Winkler.
\newblock Gibbs measures and dismantlable graphs.
\newblock {\em J. Combin. Theory Ser. B}, 78(1):141--166, 2000.

\bibitem{BW04}
Graham~R. Brightwell and Peter Winkler.
\newblock Graph homomorphisms and long range action.
\newblock In {\em Graphs, morphisms and statistical physics}, volume~63 of {\em
  DIMACS Ser. Discrete Math. Theoret. Comput. Sci.}, pages 29--47. Amer. Math.
  Soc., Providence, RI, 2004.

\bibitem{DocGro}
Anton Dochtermann.
\newblock {Homotopy groups of Hom complexes of graphs}.
\newblock arXiv:math.CO/07052620.

\bibitem{GR01}
Chris Godsil and Gordon Royle.
\newblock {\em Algebraic graph theory}, volume 207 of {\em Graduate Texts in
  Mathematics}.
\newblock Springer-Verlag, New York, 2001.

\bibitem{HN90}
Pavol Hell and Jaroslav Ne{\v{s}}et{\v{r}}il.
\newblock On the complexity of {$H$}-coloring.
\newblock {\em J. Combin. Theory Ser. B}, 48(1):92--110, 1990.

\bibitem{HN04}
Pavol Hell and Jaroslav Ne{\v{s}}et{\v{r}}il.
\newblock {\em Graphs and homomorphisms}, volume~28 of {\em Oxford Lecture
  Series in Mathematics and its Applications}.
\newblock Oxford University Press, Oxford, 2004.

\bibitem{HOV99}
Mark Hovey.
\newblock {\em Model categories}, volume~63 of {\em Mathematical Surveys and
  Monographs}.
\newblock American Mathematical Society, Providence, RI, 1999.

\bibitem{Kcol}
Dmitry~N. Kozlov.
\newblock Collapsing along monotone poset maps.
\newblock {\em Int. J. Math. Math. Sci.}, pages Art. ID 79858, 8, 2006.

\bibitem{Kozsimp}
Dmitry~N. Kozlov.
\newblock Simple homotopy types of {H}om-complexes, neighborhood complexes,
  {L}ov\'asz complexes, and atom crosscut complexes.
\newblock {\em Topology Appl.}, 153(14):2445--2454, 2006.

\bibitem{Kfold}
Dmitry~N. Kozlov.
\newblock A simple proof for folds on both sides in complexes of graph
  homomorphisms.
\newblock {\em Proc. Amer. Math. Soc.}, 134(5):1265--1270 (electronic), 2006.

\bibitem{Kchr}
Dmitry~N. Kozlov.
\newblock Chromatic numbers, morphism complexes, and {S}tiefel-{W}hitney
  characteristic classes.
\newblock In {\em Geometric combinatorics}, volume~13 of {\em IAS/Park City
  Math. Ser.}, pages 249--315. Amer. Math. Soc., Providence, RI, 2007.

\bibitem{Lov78}
L.~Lov{\'a}sz.
\newblock Kneser's conjecture, chromatic number, and homotopy.
\newblock {\em J. Combin. Theory Ser. A}, 25(3):319--324, 1978.

\bibitem{Mac98}
Saunders Mac~Lane.
\newblock {\em Categories for the working mathematician}, volume~5 of {\em
  Graduate Texts in Mathematics}.
\newblock Springer-Verlag, New York, second edition, 1998.

\bibitem{QUI67}
Daniel~G. Quillen.
\newblock {\em Homotopical algebra}.
\newblock Lecture Notes in Mathematics, No. 43. Springer-Verlag, Berlin, 1967.

\end{thebibliography}

\end{document}